\magnification=1200
\def\qed{\unskip\kern 6pt\penalty 500\raise -2pt\hbox
{\vrule\vbox to 10pt{\hrule width 4pt\vfill\hrule}\vrule}}
\centerline{STRUCTURE AND $f$-DEPENDENCE OF THE A.C.I.M.}
\centerline{FOR A UNIMODAL MAP $f$ OF MISIUREWICZ TYPE.}
\bigskip
\centerline{by David Ruelle\footnote{*}{Mathematics Dept., Rutgers University, and IHES, 91440 Bures sur Yvette, France.\break $<$ruelle@ihes.fr$>$}.}
\bigskip\bigskip\noindent
	{\leftskip=2cm\rightskip=2cm\sl {\bf Abstract.}  By using a suitable Banach space on which we let the transfer operator act, we make a detailed study of the ergodic theory of a unimodal map $f$ of the interval in the Misiurewicz case.  We show in particular that the absolutely continuous invariant measure $\rho$ can be written as the sum of 1/square root  spikes along the critical orbit, plus a continuous background.  We conclude by a discussion of the sense in which the map $f\mapsto\rho$ may be differentiable.
\par}
\vfill\eject
\noindent
{\bf 0 Introduction.}
\medskip
	This paper is part of an attempt to understand the smoothness of the map $f\mapsto\rho$ where $(M,f)$ is a differentiable dynamical system and $\rho$ an SRB measure.  [For a general introduction to the problems involved, see for instance [2], [31]].  Smoothness has been established for uniformly hyperbolic systems (see [17], [21], [14], [22], [9]).  In that case, one finds that the derivative of $\rho$ with respect to $f$ can be expressed in terms of the value at $\omega=0$ of a {\it susceptibility function} $\Psi(e^{i\omega})$ which is holomorphic when the {\it complex frequency} $\omega$ satisfies Im $\omega>0$, and meromorphic for Im $\omega>$ some negative constant.  In the absence of uniform hyperbolicity, $f\mapsto\rho$ need not be continuous.  Consider then a family $(f_\kappa)_{\kappa\in{\bf R}}$.  A theorem of H. Whitney [29] gives general conditions under which, if $\rho_\kappa$ is defined on $K\subset{\bf R}$, then $\kappa\mapsto\rho_\kappa$ extends to a differentiable function of $\kappa$ on ${\bf R}$.  Taking $\rho_\kappa$ to be an SRB measure for $f_\kappa$, this gives a reasonable meaning to the differentiability of $\kappa\mapsto\rho_\kappa$ on $K$ (as proposed in [24], see [20], [11] for a different application of Whitney's theorem), even though we start with a noncontinuous function $\kappa\mapsto\rho_\kappa$ on ${\bf R}$.\medskip
	Using Whitney's theorem to study SRB states as proposed above is a delicate matter.  A simple situation that one may try to analyze is when $(M,f)$ is a unimodal map of the interval and $\rho$ an absolutely continuous invariant measure (a.c.i.m.).  [From the vast literature on this subject, let us mention [12], [13], [6], [7], [8], [28]].  A preliminary study of the Markovian case ({\it i.e.}, when the critical orbit is finite, see [23], [16]) shows that the susceptibility function $\Psi(\lambda)$ has poles for $|\lambda|<0$, but is holomorphic at $\lambda=1$.  This study suggests that in non-Markovian situations $\Psi$ may have a natural boundary separating $\lambda=0$ (around which $\Psi$ has a natural expansion) and $\lambda=1$ (corresponding to $\omega=0$).  Misiurewicz [19] has studied a class of unimodal maps where the critical orbit stays away from the critical point, and he has proved the existence of an a.c.i.m. $\rho$ for this class.  This seems a good situation where one could study the dependence of $\rho$ on $f$, as pointed out to the author by L.-S. Young.
\medskip
	A desirable starting point to study the dependence of the a.c.i.m. $\rho$ on $f$ is to have an operator ${\cal L}$ on a Banach space ${\cal A}$ such that ${\cal L}\rho=\rho$, and $1$ is a simple isolated eigenvalue of ${\cal L}$.  The main content of the present paper is the construction of ${\cal A}$ and ${\cal L}$ with the desired properties.  Specifically we write ${\cal A}={\cal A}_1\oplus{\cal A}_2$, where ${\cal A}_2$ consists of {\it spikes}, {\it i.e.}, $1/$square root singularities at points of the critical orbit, which are known to be present in $\rho$.  We are thus able to prove that the a.c.i.m. $\rho$ is the sum of a continuous background, and of the spikes (see Theorem 9, and the Remarks 16).  Note that the construction of an operator ${\cal L}$ with a spectral gap had been achieved earlier by G. Keller and T. Nowicki [18], and by L.-S. Young [30] (our construction, in a more restricted setting, leads to stronger results).
\medskip
	We start studying the smoothness of the map $f\mapsto\rho$ by an informal discussion in Section 17.  Theorem 19 proves the differentiability along topological conjugacy classes (which are codimension 1) and relates the derivative to the value at $\lambda=1$ of a modified susceptibility function $\Psi(X,\lambda)$.  [Following an idea of Baladi and Smania [5], it is plausible that differentiability in the sense of Whitney holds in directions tangent to a conjugacy class, see below].  Transversally to topological conjugacy classes the map $f\mapsto\rho$ is continuous, but appears not to be differentiable.  While this nondifferentiability is not rigorously proved, it seems to be an unavoidable consequence of the fact that the weight of the $n$-th spike is roughly $\sim\alpha^{n/2}$ (for some $\alpha\in(0,1)$) while its speed when $f$ changes is $\sim\alpha^{-n}$.  [See Section 16(c).  In fact, for a smooth family $(f_\kappa)$ restricted to values $\kappa\in K$ such that $f_\kappa$ is in a suitable Misiurewicz class, the estimates just given for the weight and speed of the spikes suggest that $\kappa\to\rho_\kappa(A)$ for smooth $A$ is ${1\over2}$-H\"older, and nothing better, but we have not proved this].  Physically, let us remark that the spikes of high order $n$ will be drowned in noise, so that discontinuities of the derivative of $f\mapsto\rho$ will be invisible.
\medskip
	Note that the susceptibility functions $\Psi(\lambda)$, $\Psi(X,\lambda)$ to be discussed may have singularities both for large $|\lambda|$ and small $|\lambda|$.  [The latter singularities do not occur for uniformly hyperbolic systems, but show up for the unimodal maps of the interval in the Markovian case, as we have mentioned above.  A computer search of such singularities is of interest [10]].
\medskip
	A study similar to that of the present paper has been made (Baladi [3], Baladi and Smania [5]) for piecewise expanding maps of the interval.  In that case it is found that $f\mapsto\rho$ is not differentiable in general, but Baladi and Smania study the differentiability of $f\mapsto\rho$ along directions tangent to topological conjugacy classes (horizontal directions), not just for $f$ restricted to a class.  Note that our $1/$square root spikes are replaced in the piecewise expanding case by jump discontinuities.  This entails some serious differences, in particular, in the piecewise expanding case $\Psi(\lambda)$ is holomorphic for $|\lambda|<1$.
\medskip\noindent
{\bf Acknowledgments.}
\medskip
	I am very indebted to Lai-Sang Young and Viviane Baladi for their help and advice in the elaboration of the present paper.  L.-S. Young was most helpful in getting this study started, and V. Baladi in getting it finished.
\bigskip\noindent
{\bf 1 Setup.}
\medskip
	Let $I$ be a compact interval of ${\bf R}$, and $f:{\bf R}\to{\bf R}$ be real-analytic.  We assume that there is $c$ in the interior of $I$ such that $f'(c)=0$, $f'(x)>0$ for $x<c$, $f'(x)<0$ for $x>c$, and $f''(c)<0$.  Replacing $I$ by a possibly smaller interval, we assume that $I=[a,b]$ where $b=fc$, $a=f^2c$, and $a<fa$. 
\medskip
	We shall construct a {\it horseshoe} $H\subset(a,b)$, {\it i.e.}, a mixing compact hyperbolic set with a Markov partition for $f$.  Following Misiurewicz [19] we shall assume that $fa\in H$.
\medskip
	Under natural conditions to be discussed below we shall study the existence of an  a.c.i.m. $\rho(x)\,dx$ for $f$, and its dependence on $f$.
\bigskip\noindent
{\bf 2 Construction of the set $H(u_1)$.}
\medskip
	Let $u_1\in[a,b]$ and define the closed set
$$	H(u_1)=\{x\in[a,b]:f^nx\ge u_1\hbox{ for all }n\ge0\}      $$
We have thus $fH(u_1)\subset H(u_1)$.  Assuming that $H(u_1)$ is nonempty, let $v$ be its minimum element, then $H(u_1)=H(v)$.  [Since $v\in H(u_1)$ we have $v\ge u_1$, hence $H(v)\subset H(u_1)$.  If $H(u_1)$ contained an element $w\notin H(v)$ we would have $H(u_1)\ni f^kw<v$ for some $k\ge0$, in contradiction with the minimality of $v$].  Therefore we may (and shall) assume that $H(u_1)\ni u_1$.  We shall also assume
$$		a<u_1<c,fa      $$
(and $f^2u_1\not=u_1$, which will later be replaced by a stronger condition).
There is $u_2\in[a,b]$ such that $fu_2=u_1$ and, since $u_1<fa$, it follows that $u_2$ is unique and satisfies $c<u_2<b$.  We have $u_2\in H(u_1)$ [because $u_2>c>u_1$ and $fu_2\in H(u_1)$] and if $x\in H(u_1)$ then $x\le u_2$ [because $x>u_2$ implies $fx<u_1$].  Therefore, $u_2$ is the maximum element of $H(u_1)$.  Let
$$	 V_0=\{x\in[a,b]:fx>u_2\}      $$
then $u_1<V_0$ [because $x\le u_1$ implies $fx\le fu_1\in H(u_1)\le u_2$] and $V_0<u_2$ [because $x\ge u_2$ implies $fx\le fu_2=u_1<u_2$]. Thus we may write $V_0=(v_1,v_2)$, with $u_1<v_1<c<v_2<u_2$ [$u_1\not=v_1$ because $f^2u_1\not=u_1$].  We have $v_1,v_2\in H(u_1)$ [because $v_1,v_2>u_1$ and $fv_1=fv_2=u_2\in H(u_1)$].
\medskip
	Our assumptions ($H(u_1)\ni u_1$, $a<u_1<c,fa$ and $f^2u_1\not=u_1$) and definitions give thus
$$	H(u_1)\subset[u_1,v_1]\cup[v_2,u_2]      $$
$$	f[u_1,v_1]\subset[u_1,u_2]\qquad,\qquad f[v_2,u_2]=[u_1,u_2]      $$
and
$$	H(u_1)=\{x\in[u_1,u_2]:f^nx\notin V_0\hbox{ for all }n\ge0\}=fH(u_1)      $$
\indent
	Let us say that the open interval $V_\alpha\subset[u_1,u_2]$ is of order $n$ if $f^n$ maps homeomorphically $V_\alpha$ onto $(v_1,v_2)=V_0$.  We have thus
$$	H(u_1)=[u_1,u_2]\backslash\cup\hbox{ all $V_\alpha$}      $$
By induction on $n$ we shall see that 
$$	[u_1,u_2]\backslash\cup\hbox{ the $V_\alpha$ of order $\le n$}      $$
is composed of disjoint closed intervals $J$, such that $f^nJ\subset[u_1,v_1]$ or $[v_2,u_2]$ when $n>0$, and the endpoints of $f^nJ$ are $u_1,u_2,v_1,v_2$ or an image of these points by $f^k$ with $k\le n$.  Assume that the induction assumption holds for $n$ (the case of $n=0$ is trivial) and let $J$ be as indicated.  Since $f^nJ\subset[u_1,v_1]$ or $[v_2,u_2]$, $f^{n+1}$ is monotone on $J$, and the endpoints of $J$ are mapped by $f^{n+1}$ outside of $V_0$ [because $u_1,u_2,v_1,v_2$ and their images by $f^\ell$ are in $H(u_1)$, hence $\notin (v_1,v_2)$].  The interval $V_0$ is thus either inside of $f^{n+1}J$ or disjoint from $f^{n+1}J$.  Each $V_\alpha$ of order $n+1$ thus obtained is disjoint from other $V_\alpha$ of order $\le n+1$, and the closed intervals $\tilde J$ in $[u_1,u_2]\backslash\cup\hbox{ the $V_\alpha$ of order $\le n+1$}$, are such that the endpoints of $f^{n+1}\tilde J$ are $u_1,u_2,v_1,v_2$ or an image of these points by $f^k$ with $k\le n+1$, in agreement with our induction assumption.
\medskip
	We assume now that, for some $N\ge0$, we have $f^{N+1}u_1=u_1$ (take $N$ smallest with this property), and we assume also that $(f^{N+1})'(u_1)>0$.  [$N=0,1$ cannot occur, in particular $f^2u_1\not=u_1$.  Thus $N\ge2$, with $f^Nu_1=u_2$, $f^{N-1}u_1\in\{v_1,v_2\}$.  Furthermore, $(f^{N-1})'(u_1)<0$ if $f^{N-1}u_1=v_1$, and $(f^{N-1})'(u_1)>0$ if $f^{N-1}u_1=v_2$, {\it i.e.}, $f^{N-1}(u_1+)=v_1-$ or $v_2+$].
\medskip
	Using the above assumption we now show that none of the intervals $J$ in
$$	[u_1,u_2]\backslash\cup\hbox{ the $V_\alpha$ of order $\le n$}      $$
is reduced to a point.  We proceed by induction on $n$, assuming that $f^nJ=[f^nx_1,f^nx_2]$, where $f^nx_1<f^nx_2$ and $f^nx_1$ is of the form $v_2,u_1$ or $f^\ell u_1$ with $(f^\ell)'(u_1)>0$ while $f^nx_2$ is of the form $v_1,u_2$ or $f^\ell u_2$ with $(f^\ell)'(u_2)>0$.  Therefore the lower limit of $f^{n+1}J$ is of the form $f^mu_1$ with $(f^m)'(u_1)>0$ while the upper limit is of the form $f^mu_2$ with $(f^m)'(u_2)>0$.  If
$$	f^{n+1}J\supset (v_1,v_2)      $$
so that a new $V_\alpha$ of order $n+1$ is created, the set $f^{n+1}J\backslash(v_1,v_2)$ consists of two closed intervals, and one of them can be reduced to a point only if $f^mu_1=v_1$ with $(f^m)'(u_1)>0$ or if $f^mu_2=v_2$ with $(f^m)'(u_2)>0$.  So, either $f^{m+2}u_1=u_1$ with $(f^{m+2})'(u_1)<0$, or $f^{m+1}u_2=u_2$ with $(f^{m+1})'(u_2)<0$ hence $f^{m+1}u_1=u_1$ with $(f^{m+1})'(u_1)<0$, in contradiction with our assumption that $(f^{N+1})'(u_1)>0$.
\bigskip\noindent
{\bf 3 Consequences.}
\medskip
	(No isolated points)
\medskip
	$H(u_1)$ is obtained from $[u_1,u_2]$ by taking away successively intervals $V_\alpha$ of increasing order.  A given $x\in H(u_1)$ will, at each step, belong to some small closed interval $J$, and the endpoints of $J$ will not be removed in later steps, so that $x$ cannot be an isolated point: {\sl $H(u_1)$ has no isolated points}.
\medskip
	(Markov property)
\medskip
	Our assumption $f^{N+1}u_1=u_1$ implies that, for $n=1,\ldots,N-1$, the point $f^nu_1$ is one of the endpoints of an interval $V_\alpha$ of order $N-1-n$, which we call $V_{N-1-n}$.  These open intervals $V_k$ are disjoint, and their complement in $[u_1,u_2]$ consists of $N$ intervals $U_1,\ldots,U_N$.  Each $U_i$ is closed, nonempty, and not reduced to a point.  Furthermore, each $U_i$ (for $i=1,\ldots,N$) is mapped by $f$ homeomorphically to a union of intervals $U_j$ and $V_k$: this is what we call {\it Markov property}.
\medskip
	We impose now the following condition:
\bigskip\noindent
{\bf 4 Hyperbolicity}.
\medskip
	{\sl There are constants $A>0,\alpha\in(0,1)$ such that if $x,fx,\ldots,f^{n-1}x\in[u_1,v_1]\cup[v_2,u_2]$, then 
$$	\big|{d\over dx}f^nx\big|^{-1}<A\alpha^n      $$}
\indent
	We label the intervals $U_1,\ldots,U_N$ from left to right, so that $u_1$ is the lower endpoint of $U_1$, and $u_2$ the upper endpoint of $U_N$.  Define also an oriented graph with vertices $U_j$ and edges $U_j\to U_k$ when $fU_j\supset U_k$.  Write $U_{j_0}\buildrel\ell\over\Longrightarrow U_{j_\ell}$ if $U_{j_0}\to U_{j_1}\to\cdots\to U_{j_\ell}$, and $U_j\Longrightarrow U_k$ if $U_j\buildrel\ell\over\Longrightarrow U_k$ for some $\ell>0$.
\bigskip\noindent
{\bf 5 Lemma} (mixing).
\medskip
	{\sl (a) For each $U_j$ there is $r\ge0$ such that $U_j\buildrel{r+3}\over\Longrightarrow U_1$.
\medskip
	(b) If there is $s>0$ such that $U_1\buildrel s\over\Longrightarrow U_1$ and $U_1\buildrel s\over\Longrightarrow U_N$, then $U_1\buildrel s\over\Longrightarrow U_k$ for $k=1,\ldots,N$.
\medskip
   (c) If there is $s>0$ such that $U_j\buildrel s\over\Longrightarrow U_k$ for all $U_j,U_k\in\{U_j:U_1\Longrightarrow U_j\Longrightarrow U_1\}$, then $U_j\buildrel s\over\Longrightarrow U_k$ for all $U_j,U_k\in\{U_1\ldots,U_N\}$, and we say that $H(u_1)$ is {\rm mixing}.
\medskip
   (d) In particular if $N+1$ is a prime, then $H(u_1)$ is mixing.
\medskip
	(e) Let $u_1<\tilde u_1<c,fa$, and suppose that $f^{\tilde N+1}\tilde u_1=\tilde u_1$, $(f^{\tilde N+1})'(u_1)>0$.  Then if $H(u_1)$ is mixing, so is $H(\tilde u_1)$.}
\medskip
	(a) The interval  $U_j$ is contained in either $[u_1,v_1]$ or $[v_2,u_2]$.  Let the same hold for the successive images up to $f^rU_j$, but $f^{r+1}U_j\ni c$ [hyperbolicity and the fact that $U_j$ is not reduced to a point imply that $r$ is finite].  Then $U_j\buildrel{r+1}\over\Longrightarrow U_k$ with $U_k\ni v_1$ or $v_2$, hence $U_k\buildrel{2}\over\Longrightarrow U_1$ and $U_j\buildrel{r+3}\over\Longrightarrow U_1$.
\medskip
	(b) The $U_j$ such that $U_1\buildrel s\over\Longrightarrow U_j$ form a set of consecutive intervals and, since this set contains $U_1$ and $U_N$ by assumption, it contains all $U_j$ for $j=1,\ldots,N$.
\medskip
   (c) By assumption, $U_1\buildrel s\over\Longrightarrow U_1$ and $U_1\buildrel s\over\Longrightarrow U_N$, so that $U_1\buildrel s\over\Longrightarrow U_k$ for $k=1,\ldots,N$ by (b).  Therefore, $\{U_j:U_1\Longrightarrow U_j\Longrightarrow U_1\}=\{U_1,\ldots,U_N\}$ by (a), and thus $U_j\buildrel s\over\Longrightarrow U_k$ for all $U_j,U_k\in\{U_1\ldots,U_N\}$.
\medskip
	(d) The {\it transitive} set $\{U_j:U_1\Longrightarrow U_j\Longrightarrow U_1\}$ decomposes into $n$ disjoint subsets $S_0,\ldots,S_{n-1}$ such that $S_0\buildrel 1\over\Longrightarrow S_1\buildrel 1\over\Longrightarrow\cdots\buildrel 1\over\Longrightarrow S_{n-1}\buildrel 1\over\Longrightarrow S_0$ and there is $s>0$ such that $U_j\buildrel{sn}\over\Longrightarrow U_k$ for all $U_j,U_k\in S_m$, where $m=0,\ldots,n-1$.  We may suppose that $U_1\in S_0$, and therefore if $U_{(k)}$ denotes the interval containing $f^ku_1$ we have $U_{(k)}\in S_{(k)}$ where $(k)=k(\hbox{mod }n)$.  Therefore $N+1$ is a multiple of $n$, where $n\le N<N+1$.  In particular, if $N+1$ is prime, then $n=1$, and $U_j\buildrel s\over\Longrightarrow U_k$ for all $U_j,U_k\in\{U_j:U_1\Longrightarrow U_j\Longrightarrow U_1\}$, so that (c) can be applied.
\medskip
	(e) Since $H(\tilde u_1)$ is a compact subset of $H(u_1)$, without isolated points, the fact that $H(u_1)$ is mixing implies that $H(\tilde u_1)$ is mixing.\qed
\bigskip\noindent
{\bf 6 Horseshoes.}
\medskip
	Note that we have
$$	H(u_1)=\{x\in[u_1,u_2]:f^nx\notin V_0\hbox{ for all }n\ge0\}
	=\cap_{n\ge0}f^{-n}([u_1,u_2]\backslash V_0)      $$
The sets $U_i\cap H(u_1)$ form a {\it Markov partition} of $H(u_1)$, {\it i.e.}, $f(U_i\cap H(u_1))$ is a finite union of sets $U_j\cap H(u_1)$.
\medskip
	A set $H=H(u_1)$ as constructed in Section 2, with the hyperbolicity and mixing conditions will be called a {\it horseshoe}.  A horseshoe is thus a mixing hyperbolic set with a Markov partition.
\medskip
	Remember that the open interval $V_\alpha\subset[u_1,u_2]$ is of order $n$ if $f^n$ maps $V_\alpha$ homeomorphically onto $V_0=(v_1,v_2)$, and let $|V_\alpha|$ be the length of $V_\alpha$.  Hyperbolicity has the following consequence.
\bigskip\noindent
{\bf 7 Lemma} (a consequence of hyperbolicity).
\medskip
	{\sl There are constants $B>0$, $\beta\in(0,1)$ such that
$$	\sum_{\alpha:{\rm order}\,V_\alpha=n}|V_\alpha|\le B\beta^n      $$}
\indent
	It suffices to prove that
$$	\hbox{Lebesgue meas. }([u_1,u_2]\backslash\cup\hbox{the $V_\alpha$ of order}
\le n)\le G\beta^n      $$
[incidentally, this shows that $H(u_1)$ has Lebesgue measure 0].
\medskip
	Let $J$ denote one of the closed intervals in
$$	[u_1,u_2]\backslash\cup\hbox{the $V_\alpha$ of order}\le n      $$
and suppose that $J$ is one of the two intervals adjacent to a given $V_\alpha$ of order $n$.  There is $n'>n$ such that $J$ contains no interval $V$ of order $<n'$, but $J\supset V_{\alpha'}$ of order $n'$.  We write $J=J_{nn'}(V_\alpha,V_{\alpha'})$ and note that $J$ is entirely determined by $V_\alpha$ and $V_{\alpha'}$ (of orders $n,n'$ respectively).  The intervals in
$$	[u_1,u_2]\backslash\cup\hbox{the $V_\alpha$ of order}\le n      $$
are all the $J_{n_1n_2}$ with $n_1\le n$ and $n_2>n$.  There is a graph $\Gamma$ with vertices $V_\alpha$ and oriented edges $J_{nn'}(V_\alpha,V_{\alpha'})$ such that for each $V_\alpha$ of order $n$ two edges $J_{nn_1}(V_\alpha,V_{\alpha_1})$ come out of $V_\alpha$ and, if $n>0$, one edge $J_{n_0n}(V_{\alpha_0},V_\alpha)$ goes in.  The graph $\Gamma$ is a tree, rooted at $V_0$.
\medskip
	We want to show that
$$	\sum_{n_1\le n,n_2>n}\sum_{\alpha_1\alpha_2}
	|J_{n_1n_2}(V_{\alpha_1},V_{\alpha_2})|\le G\beta^n      $$
In order to do this we shall introduce intervals $\tilde J_{n_1n_2}^n(V_{\alpha_1},V_{\alpha_2})\supset J_{n_1n_2}(V_{\alpha_1},V_{\alpha_2})$ such that, for fixed $n$, the $\tilde J_{n_1n_2}^n(V_{\alpha_1},V_{\alpha_2})$ are disjoint, and we shall find $\theta\in(0,1)$ and an integer $N>0$ such that
$$	\sum_{n_1\le n,n_2>n}
	|\tilde J_{n_1+2N,n_2+2N}^{n+2N}(V_{\alpha'_1},V_{\alpha'_2})|
	\le\theta\sum_{n_1\le n,n_2>n}
	|\tilde J_{n_1n_2}^n(V_{\alpha_1},V_{\alpha_2})|      $$
(where sums over $\alpha'_1,\alpha'_2$ and $\alpha_1,\alpha_2$ are implied).  In fact, we shall prove that
$$	\textstyle{\sum^*}|\tilde J_{n'_1,n'_2}^{n+2N}(V_{\alpha'_1},V_{\alpha'_2})|
	\le\theta|\tilde J_{n_1n_2}^n(V_{\alpha_1},V_{\alpha_2})|\eqno{(*)}   $$
for fixed $\tilde J_{n_1n_2}^n(V_{\alpha_1},V_{\alpha_2})$ such that $n_1\le n,n_2>n$, where the sum $\sum^*$ extends over all $\tilde J_{n'_1,n'_2}^{n+2N}(V_{\alpha'_1},V_{\alpha'_2})$ such that $J_{n'_1,n'_2}(V_{\alpha'_1},V_{\alpha'_2})$ is above $J_{n_1n_2}^n(V_{\alpha_1},V_{\alpha_2})$ in the tree $\Gamma$, and that $n'_1\le n+2N$, $n'_2>n+2N$.  [This means that $\sum^*$ extends over $\tilde J^{n+2N}$ corresponding to the closed intervals $J^*$ of 
$$	[u_1,u_2]\backslash\cup\hbox{the $V_{\alpha'}$ of order $\le n+2N$}      $$
such that $J^*\subset J_{n_1n_2}(V_{\alpha_1},V_{\alpha_2})$].
\medskip
	Note that $J_{n_1n_2}(V_{\alpha_1},V_{\alpha_2})\supset V_{\alpha_2}$ and that for some constant $K_1$ independent of $n_1,n_2$ we may write $|J_{n_1n_2}(V_{\alpha_1},V_{\alpha_2})|\le K_1|V_{\alpha_2}|$ [otherwise $J_{n_1n_2}(V_{\alpha_1},V_{\alpha_2})$would contain a $V_\alpha$ of order $<n_2$].  We can also compare $|V_{\alpha_1}|$ and $|V_{\alpha_2}|$ because $f^{n_1}V_{\alpha_1}=f^{n_2}V_{\alpha_2}=V_0$: using hyperbolicity and the smoothness of $f$ we find a constant $K_2$ such that $|V_{\alpha_2}|\le K_2\alpha^{n_2-n_1}|V_{\alpha_1}|$.  Thus
$$	|J_{n_1n_2}(V_{\alpha_1},V_{\alpha_2})|
	\le K_1K_2\alpha^{n_2-n_1}|V_{\alpha_1}|
	\le\alpha^{n_2-n_1-N}{1\over3}|V_{\alpha_1}|      $$
for suitable $N$.  We also assume that $2\alpha^N<1$.
\medskip
	If $n_2-n_1<2N$ we define $\tilde J_{n_1n_2}^n(V_{\alpha_1},V_{\alpha_2})=J_{n_1n_2}(V_{\alpha_1},V_{\alpha_2})$.  If $n_2-n_1\ge 2N$ we define $\tilde J_{n_1n_2}^n(V_{\alpha_1},V_{\alpha_2})$ as the union of $J_{n_1n_2}(V_{\alpha_1},V_{\alpha_2})$ and an adjacent subinterval $\tilde V\subset V_{\alpha_1}$ such that $|\tilde V|=\alpha^{{1\over2}(n-n_1)}{1\over3}|V_{\alpha_1}|$ and therefore (since $n<n_2$)
$$	|\tilde V|>\alpha^{{1\over2}(n_2-n_1)}{1\over3}|V_{\alpha_1}|
	>\alpha^{n_2-n_1-N}{1\over3}|V_{\alpha_1}|
	\ge|J_{n_1n_2}(V_{\alpha_1},V_{\alpha_2})|      $$
\indent
	If $n+2N<n_2$, there is only one term in the left-hand side of ($*$), and this term is $\tilde J_{n_1n_2}^{n+2N}(V_{\alpha_1},V_{\alpha_2})$, so that
$$	\Big|{\tilde J_{n_1n_2}^{n+2N}(V_{\alpha_1},V_{\alpha_2})
	\over\tilde J_{n_1n_2}^n(V_{\alpha_1},V_{\alpha_2})}\Big|\le
	{\alpha^{{1\over2}(n-n_1+2N)}{1\over3}|V_{\alpha_1}|
	+\alpha^{n_2-n_1-N}{1\over3}|V_{\alpha_1}|
	\over\alpha^{{1\over2}(n-n_1)}{1\over3}|V_{\alpha_1}|}      $$
$$	=\alpha^N+\alpha^{n_2-{1\over2}n_1-{1\over2}n-N}
	\le\alpha^N+\alpha^{n_2-n-N}\le2\alpha^N      $$
If $n+2N\ge n_2$ there are several terms in the left-hand side of ($*$), obtained from the interval $J_{n_1n_2}(V_{\alpha_1},V_{\alpha_2})$ from which at least a subinterval of length ${1\over3}|V_{\alpha_2}|$ has been taken out.  Therefore
$$	\textstyle{\sum^*}
	\le|J_{n_1n_2}(V_{\alpha_1},V_{\alpha_2})|-{1\over3}|V_{\alpha_2}|      $$
and
$$	{\textstyle{\sum^*}\over|\tilde J_{n_1n_2}^n(V_{\alpha_1},V_{\alpha_2})|}
	\le1-{{1\over3}|V_{\alpha_2}|\over|J_{n_1n_2}(V_{\alpha_1},V_{\alpha_2})|}
	\le1-{{1\over3}|V_{\alpha_2}|\over K_1|V_{\alpha_2}|}\le1-{1\over3K_1}      $$
We have thus proved ($*$) with $\theta=\max(2\alpha^N,1-1/3K_1)$, and the lemma follows, with $\beta^N=\theta$.\qed
\bigskip\noindent
{\bf 8 Remark} (the set $\tilde H$).
\medskip
	Starting from the horseshoe $H=H(u_1)$ we can, by increasing $u_1$ to $\tilde u_1$ such that $\tilde u_1<c,fa$, obtain a set $\tilde H=H(\tilde u_1)\subset H$ such that $\tilde u_1\in\tilde H$ and the distance of $\tilde H$ to $\{u_1,u_2,v_1,v_2\}$ is $\ge\epsilon>0$.  [In fact, using our hyperbolicity assumption we can arrange that there is $\tilde N$ such that $f^{\tilde N+1}\tilde u_1=\tilde u_1$, $(f^{\tilde N+1})'(\tilde u_1)>0$.  In that case $\tilde H$ is mixing (Lemma 5(e)) and therefore again a horseshoe].
\bigskip\noindent
{\bf 9 Theorem.}
\medskip
	{\sl Let $H=H(u_1)$ be a horseshoe, suppose that $fa=f^2b\in H$, and that $\{f^nb:n\ge0\}$ has a distance $\ge\epsilon>0$ from $\{u_1,u_2,v_1,v_2\}$.  Then $f$ has a unique a.c.i.m. $\rho(x)\,dx$.  Furthermore
$$	\rho(x)=\phi(x)+\sum_{n=0}^\infty C_n\psi_n(x)      $$
The function $\phi$ is continuous on $[a,b]$, with $\phi(a)=\phi(b)=0$.  For $n\ge0$ we shall choose $w_n\in\{u_1,u_2,v_1,v_2\}$ with $(w_n-c)(c-f^nb)<0$ and let $\theta_n$ be the characteristic function of $\{x:(w_n-x)(x-f^nb)>0\}$.  Then, the above constants $C_n$ and {\rm spikes} $\psi_n$ are defined by
$$	C_n=\phi(c)|{1\over2}f''(c)\prod_{k=0}^{n-1}f'(f^kb)|^{-1/2}      $$
$$	\psi_n(x)={w_n-x\over w_n-f^nb}\cdot|x-f^nb|^{-1/2}\,\theta_n(x)      $$}
\indent
	[The condition that $\{f^nb:n\ge0\}$ has distance $\ge\epsilon$ from $\{u_1,u_2,v_1,v_2\}$ is achieved, according to Remark 8, by taking $\epsilon\le|u_1-a|,|u_2-b|$, and $f^2b\in\tilde H$.  Note also that $\psi_n(c)=0$, so that $\phi(c)=\rho(c)$.  Other choices of $\psi_n$ can be useful, with the same singularity at $f^nb$, but greater smoothness at $w_n$ and/or satisfying $\int dx\,\psi_n(x)=0$].
\bigskip\noindent
{\bf 10 Analysis.}
\medskip
	We analyze the problem before starting the proof.  Near $c$ we have
$$	y=fx=b-A(x-c)^2+\hbox{h.o.}      $$
with $A=-f''(c)/2>0$, hence $x-c=\pm\big((b-y)/A\big)^{1/2}+O(b-y)$.  Therefore, writing $U=\rho(c)/\sqrt A$, the density of the image $f(\rho(x)dx)$ by $f$ of $\rho(x)dx$ has, near $b$, a singularity
$$	{U\over\sqrt{(b-x)}}+O(\sqrt{b-x})      $$
and, near $a$, a singularity
$${	U\over\sqrt{-f'(b)(x-a)}}+O(\sqrt{x-a})      $$
To deal with the general case of the singularity at $f^nb$, define $s_n=-\hbox{sgn}\prod_{k=0}^{n-1}f'(f^kb)$, so that
$$	\prod_{k=0}^{n-1}f'(f^kb)=-s_nU^2C_n^{-2}      $$
The density of $f(\rho(x)dx)$ has then, near $f^nb$, a singularity given when $s_n(x-f^nb)>0$ by
$$	{U\over\sqrt{(\prod_{k=0}^{n-1}|f'(f^kb)|)|x-f^nb|}}+O(\sqrt{|x-f^nb|})    $$
$$	={U\over\sqrt{-(x-f^nb)\prod_{k=0}^{n-1}f'(f^kb)}}
	+O(\sqrt{|x-f^nb|})    $$
$$	={C_n\over\sqrt{s_n(x-f^nb})}+O(\sqrt{s_n(x-f^nb}))      $$
and by $0$ when $s_n(x-f^nb)<0$.
\medskip
	We let now $w_0=u_2$ and, for $n\ge0$, define $w_{n+1}\in\{u_1,u_2,v_1,v_2\}$ inductively by:
$$	(w_{n+1}-c)(f^{n+1}b-c)>0\qquad,\qquad(w_{n+1}-f^{n+1}b)(fw_n-f^{n+1}b)>0      $$
We have thus $w_0=u_2,w_1=u_1$, and in general
$$	w_n\in\{u_1,u_2,v_1,v_2\}\quad,\quad(w_n-c)(f^nb-c)>0\quad,\quad 
s_n(w_n-f^nb)>0\quad,\quad|w_n-f^nb|\ge\epsilon   $$
\medskip
	The above considerations show that the singularity expected near $f^nb$ for the density of $f(\rho(x)dx)$ is also represented by
$$	(1-{x-f^nb\over w_n-f^nb})\cdot{C_n\over\sqrt{s_n(x-f^nb)}}
	\,\theta_n(x)      $$
$$	=C_n\,{w_n-x\over w_n-f^nb}\,|x-f^nb|^{-1/2}\,\theta_n(x)
	=C_n\psi_n(x)      $$
in agreement with the claim of the theorem.
\bigskip\noindent
{\bf 11 Lemma.}
\medskip
	{\sl Write
$$	f(\psi_n(x)dx)=\tilde\psi_{n+1}(x)dx\qquad,\qquad
	\tilde\psi_{n+1}=|f'(f^nb)|^{-1/2}\psi_{n+1}+\chi_n      $$
Then, for $n\ge0$, the $\chi_n$ are continuous of bounded variation on $[a,b]$, with $\chi_n(a)=\chi_n(b)=0$, and the ${\rm Var}\,\chi_n=\int_a^b|d\chi_n/dx|dx$ are bounded uniformly with respect to $n$.  Furthermore, if $n\ge1$ and $V_\alpha\subset{\rm supp}\chi_n$, then $\chi_n|V_\alpha$ extends to a holomorphic function $\chi_{n\alpha}$ in a complex neighborhood $D_\alpha$ of the closure of $V_\alpha$ in ${\bf R}$ (further specified in Section 12), with the $|\chi_{n\alpha}|$ uniformly bounded.}
\medskip
	The case $n=0$ can be handled by inspection, and we shall assume $n\ge1$.  We let
$$	I_n=\Big\{\matrix{(fa,b)&\quad{\rm if}\quad&f^nb\in[a,c)\cr
	(a,b)&\quad{\rm if}\quad&f^nb\in(c,b)\cr}      $$
And define $f_n^{-1}:I_n\mapsto(a,b)$ to be the inverse of $f$ restricted respectively to $(a,c)$ or $(c,b)$ in the two cases above.  We have then
$$	\tilde\psi_{n+1}(x)={\psi_n(f_n^{-1}x)\over|f'(f_n^{-1}x)|}      $$
\indent
	Since $n\ge1$, the region of interest $f{\rm supp}\psi_n\cup{\rm supp}\psi_{n+1}$ is $\subset[u_1,u_2]\subset(a,b)$, and we have
$$	f_n^{-1}x-f^nb=(x-f^{n+1}b)A_n(x)      $$
where $A_n$ is real analytic and $A_n(f^{n+1}b)=(f'(f^nb))^{-1}$.  Therefore we may write
$$	{1\over f_n^{-1}x-f^nb}
	={f'(f^nb)\over x-f^{n+1}b}(1+(x-f^{n+1}b)\tilde A_n(x))      $$
$$	{1\over f'(f_n^{-1}x)}={1\over f'(f^nb)}(1+(x-f^{n+1}b))\tilde B_n(x)      $$
$$	{w_n-f_n^{-1}x\over w_n-f^nb}=1+(x-f^{n+1}b)\tilde C_n(x)      $$
and since
$$	\psi_n(f_n^{-1}x)=\theta_n(f_n^{-1}x)\Big|{w_n-f_n^{-1}x\over w_n-f^nb}	\Big|\cdot|f_n^{-1}x-f^nb|^{-1/2}      $$
we find
$$	\tilde\psi_{n+1}(x)={\theta_n(f_n^{-1}x)|f'(f^nb)|^{-1/2}\over\sqrt{|x-f^{n+1}b|}}
	\big(1+(x-f^{n+1}b)\tilde D_n(x)\big)      $$
with $\tilde D_n$ real analytic.  Note that $\tilde\psi_{n+1}$ and $|f'(f^nb)|^{-1/2}\psi_{n+1}$ have the same singularity at $f^{n+1}b$.  It follows readily that $\tilde\psi_{n+1}-|f'(f^nb)|^{-1/2}\psi_{n+1}$ is a continuous function $\chi_n$ vanishing at the endpoints of its support, and bounded uniformly with respect to $n$.  It is easy to see that ${\rm Var}\,\chi_n$ is bounded uniformly in $n$.  The extension of $\chi_n|V_\alpha$ to holomorphic $\chi_{n\alpha}$ in $D_\alpha$ is also handled readily (see Section 12 for the description of the $D_\alpha$).\qed
\bigskip\noindent
{\bf 12 The operator ${\cal L}$ and the space ${\cal A}$.}
\medskip
	We have $f(\rho(x)\,dx)=({\cal L}_{(1)}\rho)(x)\,dx$, where the transfer operator ${\cal L}_{(1)}$ on $L^1(a,b)$ is defined by 
$$	{\cal L}_{(1)}\rho=\sum_\pm{\rho\circ f_\pm^{-1}\over|f'\circ f_\pm^{-1}|}      $$
and we have denoted by
$$	f_-^{-1}:[fa,b]\mapsto[a,c]\qquad{\rm and}\qquad f_+^{-1}[a,b]\mapsto[c,b]       $$
the branches of the inverse of $f$.  The invariance of $\rho(x)\,dx$ under $f$ is thus expressed by
$$	\rho={\cal L}_{(1)}\rho      $$
We shall look for a solution of this equation in a Banach space ${\cal A}$ defined below.  Roughly speaking, ${\cal A}$ consists of functions
$$	\phi+\sum_{n=0}^\infty c_n\psi_n      $$
where the $\psi_n$ are defined in the statement of Theorem 9, and $\phi:[a,b]\to{\bf C}$ is a less singular rest  with certain analyticity properties.
\medskip
	Remember that we may write
$$	[a,b]=H\cup[a,u_1)\cup(u_2,b]\cup\hbox{the $V_\alpha$ of all orders $\ge0$}      $$
We have (see Remark 8)
$$	{\rm clos}\,[a,u_1)\subset[a,\tilde u_1)\qquad,\qquad{\rm clos}\,(u_2,b]\subset(\tilde u_2,b]
	\qquad,\qquad{\rm clos}\,V_0\subset\tilde V_0      $$
where $\tilde u_2$ and $\tilde V_0=(\tilde v_1,\tilde v_2)$, are defined for $\tilde H$ as $u_2$ and $V_0$ were defined for $H$.  It is convenient to define $V_{-1}=(u_2,b]$ and $V_{-2}=[a,u_1)$ (of order $-1$ and $-2$ respectively) so that
$$	[a,b]=H\cup\hbox{the $V_\alpha$ of all orders $\ge-2$}      $$
We also define $\tilde V_{-1}=(\tilde u_2,b],\tilde V_{-2}=[a,\tilde u_1)$.  We let now $\tilde V_\alpha$ denote the unique interval in $[a,b]\backslash\tilde H$ such that $V_\alpha\subset\tilde  V_\alpha$.  Note that the map $V_\alpha\mapsto\tilde V_\alpha$ is not injective!  
\medskip
	For each $V_\alpha$ of order $\ge0$ we may choose an open set $D_\alpha\subset{\bf C}$ such that
$$	\tilde V_\alpha\supset D_\alpha\cap{\bf R}\supset{\rm clos}\,V_\alpha      $$
and, if $fV_\beta=V_\alpha$ of order $\ge0$, $fD_\beta\supset{\rm clos}\,D_\alpha$ [we have here denoted by ${\rm clos}\,V_\alpha$ the closure of $V_\alpha$ in ${\bf R}$, and by ${\rm clos}\,D_\alpha$ the closure of $D_\alpha$ in ${\bf C}$].  Let also $R_a,R_b$ be two-sheeted Riemann surfaces, branched respectively at $a,b$, with natural projections $\pi_a,\pi_b:R_a,R_b\to{\bf C}$.  We may choose open sets $D_{-1},D_{-2}\subset{\bf C}$ such that, for $\alpha=-1,-2$,
$$	\tilde V_\alpha\supset D_\alpha\cap\{x\in{\bf R}:a\le x\le b\}\supset{\rm clos}\,V_\alpha   $$
and $f$ extends to holomorphic maps $\tilde f_{-1}:D_0\to R_b,\tilde f_{-2}:(\tilde f_{-1}D_0)\to R_a$ such that $\tilde f_{-1}D_0\supset\pi_b^{-1}{\rm clos}\,D_{-1},\tilde f_{-2}\pi_b^{-1}D_{-1}\supset\pi_a^{-1}{\rm clos}\,D_{-2}$. [We shall say that $\tilde f_{-1}$ sends $(v_1,c)$ to the {\it upper} sheet of $R_b$ and $(c,v_2)$ to the {\it lower} sheet of $R_b$; $\tilde f_{-2}$ sends the upper (lower) sheet of $R_b$ to the upper (lower) sheet of $R_a$].
\medskip
	We come now to a precise definition of the complex Banach space ${\cal A}$.  We write ${\cal A}={\cal A}_1\oplus{\cal A}_2$ where the elements of ${\cal A}_1$ are of the form $(\phi_\alpha)$ and the elements of ${\cal A}_2$ of the form $(c_n)$.  Here the index set of the $\phi_\alpha$ is the same as the index set of the intervals $V_\alpha$ (of order $\ge-2$); the index $n$ of the $c_n\in{\bf C}$ takes the values $0,1,\ldots$ [the $c_n$ should not be confused with the critical point $c$].  We assume that $\phi_\alpha$ is a holomorphic function in $D_\alpha$ when $V_\alpha$ is of order $\ge0$, while $\phi_{-1},\phi_{-2}$ are holomorphic on $\pi_b^{-1}D_{-1},\pi_a^{-1}D_{-2}$ and, for all $\alpha$, $||\phi_\alpha||=\sup_{z\in D_\alpha}|\phi_\alpha(z)|<\infty$.  
\medskip	
	[We shall later consider a function $\phi:[a,b]\to{\bf C}$ such that $\phi|V_\alpha=\phi_\alpha|V_\alpha$ when $V_\alpha$ is of order $\ge0$.  For $x\in V_{-1}$ we shall require $\phi(x)=\Delta\phi(x)=\phi_{-1}(x^+)-\phi_{-1}(x^-)$ where $x^+(x^-)$ is the preimage of $x$ by $\pi_b$ on the upper (lower) sheet of $\pi_b^{-1}D_{-1}$; for $x\in V_{-2}$ we shall require $\phi(x)=\Delta\phi_{-2}(x)=\phi_{-2}(x^+)-\phi_{-2}(x^-)$ where $x^+(x^-)$ is the preimage of $x$ by $\pi_a$ on the upper (lower) sheet of $\pi_a^{-1}D_{-2}$.  But at this point we discuss an operator ${\cal L}$ on ${\cal A}$ instead of the transfer operator ${\cal L}_{(1)}$ acting on functions $\phi+\sum_nc_n\psi_n$].
\medskip
	Let $\gamma,\delta$ be such that $1<\gamma<\beta^{-1},1<\delta<\alpha^{-1/2}$ with $\beta$ as in Lemma 7 and $\alpha$ as in the definition of hyperbolicity (Section 4).  We write
$$	||(\phi_\alpha)||_1=\sup_{n\ge-2}\gamma^n\sum_{\alpha:{\rm order}\,V_\alpha=n}
	|V_\alpha|.||\phi_\alpha||\qquad,\qquad||(c_n)||_2=\sup_{n\ge0}\delta^n|c_n|      $$
and, for $\Phi=((\phi_\alpha),(c_n))$, we let $||\Phi||=||(\phi_\alpha)||_1+||(c_n)||_2$.  We let then ${\cal A}_1,{\cal A}_2$ be the Banach spaces of sequences $(\phi_\alpha),(c_n)$ as above , such that the norms $||(\phi_\alpha)||_1,||(c_n)||_2$ are finite.  We shall define ${\cal L}$ on ${\cal A}$ such that ${\cal L}\Phi=\tilde\Phi$.  We first describe what contribution each $\phi_\alpha$ or $c_n$ gives to $\tilde\Phi$ and then we shall check that this is a consistent description of an element $\tilde\Phi$ of ${\cal A}$.
\medskip
	(i) $\displaystyle\phi_\beta\Rightarrow\hat\phi_{\beta\alpha}={\phi_\beta\over|f'|}\circ(f|D_\beta)^{-1}\qquad{\rm in}\,\,D_\alpha$ if order $\beta>0$ and $fV_\beta=V_\alpha$
\par\noindent[we have here denoted by $|f'|$ the holomorphic function $\pm f'$ such that $\pm f'>0$ for real argument, we shall use the same notation in (ii)-(vi) below].
\medskip
	(ii) $\displaystyle\phi_0\Rightarrow\Big(\hat c_0=C_0\phi_0(c)\,,\,\hat\phi_{-1}=\pm{\phi_0\over|f'|}\circ\tilde f_{-1}^{-1}-C_0\phi_0(c)(\pm{1\over2}\psi_0\circ\pi_b)\qquad{\rm in}\,\,\pi_b^{-1}D_{-1}\Big)$ where the signs $\pm$ correspond to the upper/lower sheet of $\pi_b^{-1}D_{-1}$.  We claim that $\hat\phi_{-1}$ is holomorphic in $\pi_b^{-1}D_{-1}$ as the difference of two meromorphic functions with a simple pole at the branch point $b$, with the same residue.  To see this we uniformize $\pi_b^{-1}D_{-1}$ by the map $u\mapsto b-u^2$.  We have thus to express $\displaystyle\pm{\phi_0\over|f'|}(c+x)={\phi_0\over f'}(c+x)$ in terms of $u$ where $c+x=\tilde f_{-1}^{-1}(b-u^2)$ or $u=\sqrt{b-\tilde f_{-1}(c+x)}$ which gives a meromorphic function with a simple pole $1/2\sqrt{A}u$.  Since $\pm C_0\phi_0(c)\psi_0(b-u^2)$ is meromorphic with the same simple pole, $\hat\phi_{-1}$ is holomorphic in $\pi_b^{-1}D_{-1}$.
\medskip
	(iii) $\displaystyle\phi_{-1}\Rightarrow\hat\phi_{-2}={\phi_{-1}\over|f'|}\circ\tilde f_{-2}^{-1}\qquad{\rm in}\,\,\pi_a^{-1}D_{-2}$.
\medskip
	(iv) $\displaystyle\phi_{-2}\Rightarrow\hat\phi_\alpha={\Delta\phi_{-2}\over f'}\circ f^{-1}\qquad{\rm in}\,\,D_\alpha\,\,{\rm if}\,\,f(a,u_1)\supset V_\alpha\,\,,0\,\,\hbox{otherwise}$
\par\noindent[we have written $\Delta\phi_{-2}(x)=\phi_{-2}(x^+)-\phi_{-2}(x^-)$ where $x^+(x^-)$ is the preimage of $x$ by $\pi_a$ on the upper (lower) sheet of $\pi_a^{-1}D_{-2}$].
\medskip
   (v) $\displaystyle c_0\Rightarrow\Big(\hat c_1=|f'(b)|^{-1/2}c_0\,,\,\chi_0=\pm{1\over2}c_0\big({\psi_0\over|f'|}\circ\pi_b\circ\tilde f_{-2}^{-1}-|f'(b)|^{-1/2}\psi_1\circ\pi_a\big)$ $\hbox{in $\pi_a^{-1}D_{-2}$}\Big)$ where the sign $\pm$corresponds to the upper/lower sheet of $\pi_a^{-1}D_{-2}$.
\medskip
	(vi) $\displaystyle c_n\Rightarrow\Big(\hat c_{n+1}=|f'(f^nb)|^{-1/2}c_n\,,\,
\chi_{n\alpha}=c_n\big[{\psi_n\over|f'|}\circ f_n^{-1}-|f'(f^nb)|^{-1/2}\psi_{n+1}\big]$\par
$\qquad{\rm in}\,\,D_\alpha\,\,{\rm if}\,\,V_\alpha\subset\{x:\theta_n(f_n^{-1}x)>0\}\,,0\hbox{ otherwise}\Big)$\par\noindent
if $n\ge1$.
\medskip
	We may now write
$$	\tilde\Phi=((\tilde\phi_\alpha),(\tilde c_n))      $$
where
\medskip
	$\tilde\phi_{-2}=\hat\phi_{-2}+\chi_0\qquad\hbox{(see (iii),(v))}$

	$\tilde\phi_{-1}=\hat\phi_{-1}\qquad\hbox{(see(ii))}$

	$\displaystyle\tilde\phi_\alpha=\sum_{\beta:fV_\beta=V_\alpha}\hat\phi_{\beta\alpha}+\hat\phi_\alpha+\sum_{n\ge1}\chi_{n\alpha}\hbox{ if order $\alpha\ge0$}\qquad\hbox{(see (i),(iv),(vi))}$

	$\tilde c_0=\hat c_0\qquad\hbox{(see (ii))}$

	$\tilde c_1=\hat c_1\qquad\hbox{(see (v))}$

	$\tilde c_n=\hat c_n\qquad\hbox{for $n>1$}\qquad\hbox{(see (vi))}$
\medskip
	Note that, corresponding to the decomposition ${\cal A}={\cal A}_1\oplus{\cal A}_2$, we have 
$$	{\cal L}=\pmatrix{{\cal L}_0+{\cal L}_1&{\cal L}_2\cr{\cal L}_3&{\cal L}_4\cr}      $$
where
\medskip
	${\cal L}_0(\phi_\alpha)=(\sum_{\beta:fV_\beta=V_\alpha}\hat\phi_{\beta\alpha})$

	${\cal L}_1(\phi_\alpha)=(\hat\phi_\alpha)$

	${\cal L}_2(c_n)=(\chi_0,(\sum_{n\ge1}\chi_{n\alpha})_{\alpha>-1})$
	
    ${\cal L}_3(\phi_\alpha)=(\hat c_0,(0)_{n>0})$

	${\cal L}_4(c_n)=(0,(\hat c_n)_{n>0})$
\medskip\noindent
Holomorphic functions in $D_\alpha$ are defined by (i),(iv),(vi) when order $\alpha\ge0$, and in $\pi_b^{-1}D_{-1}$, $\pi_a^{-1}D_{-2}$ by (ii),(iii),(v).  Using Lemma 7, one sees that ${\cal L}_0,{\cal L}_1$ are bounded ${\cal A}_1\to{\cal A}_1$.  Using Lemma 11, one sees that ${\cal L}_3$ is bounded ${\cal A}_2\to{\cal A}_1$.  It is also readily seen that ${\cal L}_2,{\cal L}_4$ are bounded, so that ${\cal L}:{\cal A}\to{\cal A}$ is bounded.
\bigskip\noindent
{\bf 13 Theorem} (structure of ${\cal L}$).
\medskip
	{\sl With our definitions and assumptions, the bounded operator ${\cal L}:{\cal A}\to{\cal A}$ is a compact perturbation of ${\cal L}_0\oplus{\cal L}_4$; its essential spectral radius is $\le\max(\gamma^{-1},\delta\alpha^{1/2})$.}
\medskip
	Since $fa\in\tilde H$, we may assume that $f(a,u_1)\supset V_\alpha$ implies $f(D_{-2}\backslash\hbox{negative reals})\supset {\rm clos}\,D_\alpha$.  Therefore, $\phi_{-2}\mapsto\hat\phi_\alpha|D_\alpha$ is compact.  For $N$ positive integer, define the operator ${\cal L}_{N1}$ such that
$$	{\cal L}_{N1}(\phi_\alpha)={\Delta\phi_{-2}\over f'}\circ f^{-1}\qquad
	\hbox{in $D_\alpha$ if $f(a,u_1)\supset V_\alpha$ and order $\alpha>N$ , 
	$0$ otherwise}      $$
Then ${\cal L}_1$ is a perturbation of ${\cal L}_{N1}$ by a compact operator and, using Lemma 7, we see that
$$	||{\cal L}_{N1}(\phi_\alpha)||_1\le C\sup_{n>N}\gamma^n\beta^n\to0\qquad
\hbox{when $N\to\infty$}      $$
We can write ${\cal L}_2={\cal L}_{N2}+$ finite range, where
$$	{\cal L}_{N2}(c_n)=(0,0,(\sum_{n\ge N}\chi_{n\alpha})_{\alpha\ge0})      $$
Using Lemma 11 we find a bound $||\sum_{n\ge N}\chi_{n\alpha}||\le C'\delta^N$ and, using Lemma 7, 
$$	||{\cal L}_{N2}||_{{\cal A}_2\to{\cal A}_1}\le C''\delta^N\to0\qquad
	\hbox{when $N\to\infty$}      $$
The operator ${\cal L}_3$ has one-dimensional range.  Therefore ${\cal L}_1,{\cal L}_2,{\cal L}_3$ are compact operators, and the essential spectral radius of ${\cal L}$ is the max of the essential spectral radius of ${\cal L}_0$ on ${\cal A}_1$ and ${\cal L}_4$ on ${\cal A}_2$.
\medskip
	The spectral radius of ${\cal L}_4$ is
$$	\le||{\cal L}_4^N||^{1/N}
	\le\big(\delta^N C'''\sup_{\ell\ge0}\prod_{k=0}^{N-1}|f'(f^{k+\ell}b)|^{-1/2}\big)^{1/N}
	\qquad\hbox{with limit $<\delta\alpha^{1/2}$ when $N\to\infty$}      $$
The essential spectral radius of ${\cal L}_0$ is 
$$	\le\lim_{N\to\infty}
	{\sup_{n\ge N}\gamma^n\sum_{\alpha:{\rm order}V_\alpha=n}|V_\alpha|.
	||\sum_{\beta:fV_\beta=V_\alpha}\hat\phi_{\beta\alpha}||\over 
	\sup_{n\ge N}\gamma^{n+1}\sum_{\beta:{\rm order}V_\beta=n+1}|V_\beta|.
	||\phi_\beta||}      $$
$$	\le\gamma^{-1}\lim_{{\rm order}V_\alpha\to\infty}
	{|V_\alpha|.||\sum_{\beta:fV_\beta=V_\alpha}\hat\phi_{\beta\alpha}||\over
	\sum_{\beta:fV_\beta=V_\alpha}|V_\beta|.||\phi_\beta||}=\gamma^{-1}      $$
In fact, no eigenvalue of ${\cal L}_0$ can be $>\gamma^{-1}$, so the spectral radius of ${\cal L}_0$ acting on ${\cal A}_1$ is $\le\gamma^{-1}$.  The essential spectral radius of ${\cal L}$ is thus $\le\max(\gamma^{-1},\delta\alpha^{1/2})$ .\qed
\medskip\noindent
[Note also that when $\gamma\to\beta^{-1},\delta\to1$, we have $\max(\gamma^{-1},\delta\alpha^{1/2})\to\max(\beta,\alpha^{1/2})$].
\bigskip\noindent
{\bf 14 The eigenvalue $1$ of ${\cal L}$.}
\medskip
	Let the map $\Delta:{\cal A}_1\to L^1(a,b)$ be such that $\Delta(\phi_\alpha)|(a,u_1)=\Delta\phi_{-2}$, $\Delta(\phi_\alpha)|(u_2,b)=\Delta\phi_{-1}$, and $\Delta(\phi_\alpha)|V_\beta=\phi_\beta$ if order $\beta\ge0$.  We also define $w:{\cal A}\to L^1(a,b)$ by $w((\phi_\alpha),(c_n))=\Delta(\phi_\alpha)+\sum_{n=0}^\infty c_n\psi_n$ and check readily that
$$	w{\cal L}\Phi={\cal L}_{(1)}w\Phi      $$
\indent
	If $\lambda^0\ne0$ is an eigenvalue of ${\cal L}$, and $\Phi^0=((\phi_\alpha^0),(c_n^0))$ is an eigenvector to this eigenvalue, we have $w\Phi^0\ne0$ [because $w\Phi^0=0$ implies $\phi_0^0=0$, hence $\phi_{-1}^0=0,\phi_{-2}^0=0$, and $(c_n^0)=0$; then $\Delta(\phi_\alpha^0)=0$, so $\phi_\alpha^0=0$ when order $\alpha\ge0$, {\it i.e.}, $\Phi_0=0$].  Therefore
$$	\lambda^0w\Phi^0={\cal L}_{(1)}(w\Phi^0)      $$
$$	|\lambda^0|\int_a^b|w\Phi^0|=\int_a^b|{\cal L}_{(1)}(w\Phi^0)|\le
	\int_a^b{\cal L}_{(1)}|w\Phi^0|=\int_a^b|w\Phi^0|      $$
hence $|\lambda^0|\le1$.
\medskip
	If $c_0^0=0$, then $(c_n^0)=0$, and $\lambda^0$ is thus an eigenvalue of ${\cal L}_0$ acting on ${\cal A}_1$, so that $|\lambda^0|\le\gamma^{-1}$ (see Section 13).  Therefore $|\lambda^0|>\gamma^{-1}$ implies $c_0^0\ne0,c_1^0\ne0$, hence $\Delta\phi_{-1}+c_0\psi_0\ne0$, $\Delta\phi_{-2}+c_1\psi_1\ne0$.  Note that, by analyticity, $\Delta\phi_{-2}+c_1\psi_1$ is nonzero almost everywhere in $(a,u_1)$.  The image $f(a,u_1)$ contains some (small) interval $U_{i_0}\cap f^{-1}(U_{i_1}\cap f^{-1}(U_{i_2}\ldots))$ on which the image of $\Delta\phi_{-2}+c_1\psi_1$ by ${\cal L}_{(1)}$ does not vanish, and therefore (by mixing),
$$	\int_a^b|{\cal L}_{(1)}w\Phi^0|<\int_a^b{\cal L}_{(1)}|w\Phi^0|      $$
when $w\Phi^0/|w\Phi^0|$ is not constant on $(a,b)$.  Thus either (after multiplication of $\Phi^0$ by a suitable constant $\ne0$), $w\Phi^0\ge0$, or
$$	|\lambda^0|\int_a^b|w\Phi^0|<\int_a^b|w\Phi^0|\eqno{(*)}      $$
{\it i.e.}, $|\lambda^0|<1$.  Thus $1$ is the only possible eigenvalue $\lambda^0$ with $|\lambda^0|=1$, but $1$ is an eigenvalue, otherwise the spectral radius of ${\cal L}$ would be $<1$ [contradicting the fact that $\int_a^bw{\cal L}^n\Phi=\int_a^bw\Phi>0$ when $w\Phi>0$].  ${(*)}$ also implies that if ${\cal L}\Phi^1=\Phi^1$, then $w\Phi^1$ is proportional to $w\Phi^0$, hence $\phi_0^1$ is proportional to $\phi_0^0$, hence $\Phi^1$ is proportional to $\Phi^0$.  Furthermore, the generalized eigenspace to the eigenvalue $1$ contains only the multiples of $\Phi_0$ [otherwise there would exist $\Phi^1$ such that ${\cal L}^n\Phi^1=\Phi^1+n\Phi^0$, contradicting $\int_a^bw{\cal L}\Phi^1=\int_a^bw\Phi^1$].  We have proved the first part of the following
\bigskip\noindent
{\bf 15 Proposition.}
\medskip
	{\sl (a) Apart from the simple eigenvalue $1$, the spectrum of ${\cal L}$ has radius $<1$.  The eigenvector $\Phi^0$ to the eigenvalue $1$ (after multiplication by a suitable constant $\ne0$) satisfies $w\Phi^0\ge0$.
\medskip
	(b) Write $\Phi^0=((\phi_\alpha^0),(c_n^0))$ and $\Delta(\phi_\alpha^0)=\phi^0$, then $\phi^0$ is continuous, of bounded variation, and $\phi^0(a)=\phi^0(b)=0$.}
\medskip
	The interval $[u_1,u_2]$ is divided into $N$ closed intervals $W_1,\ldots,W_N$ by the points $f^nu_1$ for $n=1,\ldots,N-1$.  The intervals $W_1,\ldots,W_N$ are ordered from left to right, by doubling the common endpoints we make the $W_j$ disjoint.  Define $\gamma^0=(\gamma_j^0)_{j=1}^N$ by $\gamma_j^0=\phi^0|W_j\in L^1(W_j)$.  Then, the equation $\Phi^0={\cal L}\Phi^0$ implies
$$	\gamma^0={\cal L}_*\gamma^0+\eta\eqno{(*)}      $$
or
$$	\gamma_j^0=\sum_k{\cal L}_{jk}\gamma_k^0+\eta_j      $$
where ${\cal L}=({\cal L}_{jk})$ is a transfer operator defined as follows.  Letting $(f^{-1})_{kj}:W_j\to W_k$ be such that $f\circ(f^{-1})_{kj}$ is the identity on $W_j$ we write
$$	{\cal L}_{jk}\gamma_k=\Big\{\matrix{{\gamma_k\circ(f^{-1})_{kj}\over|f'\circ(f^{-1})_{kj}|}&\,{\rm if}\,fW_k\supset W_j\cr 0&\,{\rm otherwise}\,\cr}      $$
[the term ${\cal L}_*\gamma^0$ in $(*)$ comes from (i) in Section 12].  We let
$$	\eta_j=\sum_{n=0}^\infty\eta_{jn}      $$
Here
$$	\eta_{j0}(x)={\Delta\phi_{-2}^0(y)\over f'(y)}      $$
if $f(a,u_1)\cap W_j$ contains more than one point, and $y\in(a,u_1),fy=x\in W_j$; we let $\eta_{j0}(x)=0$ otherwise [this term comes from (iv) in Section 12].  For $n\ge1$, we let $\eta_{jn}=C_n\chi_n|W_j$ where $\chi_n=(\psi_n/|f'|)\circ f_n^{-1}-|f'(f^nb)|^{-1/2}\psi_{n+1}$ [this term comes from (vi) in Section 12].
\medskip
	Because $fu_1$ is one of the division points between the intervals $W_j$, the function $\eta_{j0}$ is continuous on $W_j$; the $\eta_{jn}$ for $n\ge1$ are also continuous.  Furthermore, $\eta_{j0}$ and the $\eta_{jn}$ for $n\ge1$ are uniformly of bounded variation.  If ${\cal H}_j$ denotes the Banach space of continuous functions of bounded variation on $W_j$ we have thus $\eta_j\in{\cal H}_j$ for $j=1,\ldots,N$.  We shall now obtain an upper bound on the essential spectral radius of ${\cal L}_*$ acting on ${\cal H}=\oplus_1^N{\cal H}_j$ by studying $||{\cal L}_*^n-F_n||$, where $F_n$ has finite-dimensional range (we use here a simple case of an argument due to Baladi and Keller [4]).  Define
$$	W_{i_n\cdots i_0}=\{x\in W_{i_n}:fx\in W_{i_{n-1}},\ldots,f^nx\in W_{i_0}\}      $$
when $fW_{i_k}\supset W_{i_{k-1}}$ for $k=n,\ldots,1$.  For $\eta=(\eta_j)\in{\cal H}$, we let $\pi_n\eta=(\pi_{jn}\eta_j)$ where $\pi_{jn}\eta_j$ is a piecewise affine function on $W_j$ such that $(\pi_{jn}\eta_j)(x)=\eta_j(x)$ whenever $x$ is an endpoint of $W_j$ or of an interval $W_{ji_{n-1}\cdots i_0}$, and is affine between all such endpoints.  Then $F_n={\cal L}_*^n\pi_n$ has finite rank ({\it i.e.}, finite-dimensional range), and ${\cal L}_*^n-F_n={\cal L}_*^n(1-\pi_n)$ maps ${\cal H}$ to ${\cal H}$.  Let ${\rm Var}\,\gamma=\sum_1^N{\rm Var}_j\gamma_j$ where ${\rm Var}_j$ is the total variation on $W_j$.  Let also $||\cdot||_0$ denote the sup-norm and $||\cdot||=\max\{{\rm Var}\,\cdot,||\cdot||_0\}$ be the bounded variation norm.  We have
$$	{\rm Var}(\gamma-\pi_n\gamma)\le2{\rm Var}\,\gamma      $$
$$	\sum_{i_0\cdots i_n}||(\gamma-\pi_n\gamma)|W_{i_n\cdots i_0}||_0
	\le{\rm Var}\,\gamma      $$
[the second inequality follows from the first because $\gamma-\pi_n\gamma$ vanishes at the endpoints of $W_{i_n\cdots i_0}$].  Since ${\cal L}_*^n(1-\pi_n)\gamma$ vanishes at the endpoints of the $W_j$, we have
$$	||({\cal L}_*^n-F_n)\gamma||={\rm Var}(({\cal L}_*^n-F_n)\gamma)      $$
$$	={\rm Var}\sum_{i_0\cdots i_n}((\gamma-\pi_n\gamma)_{i_n}\circ\tilde f_{i_n\cdots i_0})(\tilde f'\circ\tilde f_{i_n\cdots i_0})\cdots(\tilde f'\circ\tilde f_{i_1i_0})      $$
where we have written
$$	\tilde f_{i_\ell\cdots i_0}=(f^{-1})_{i_\ell i_{\ell-1}}\circ\cdots(f^{-1})_{i_1i_0}      $$
and
$$	\tilde f'={1\over|f'|}      $$
hence
$$	||({\cal L}_*^n-F_n)\gamma||\le\sum_{i_0\cdots i_n}{\rm Var}[((\gamma-\pi_n\gamma)_{i_n}\circ\tilde f_{i_n\cdots i_0})(\tilde f'\circ\tilde f_{i_n\cdots i_0})\cdots(\tilde f'\circ\tilde f_{i_1i_0})]      $$
$$	=\sum_{i_0\cdots i_n}{\rm Var}[((\gamma-\pi_n\gamma)|W_{i_n\cdots i_0})
	\prod_{\ell=0}^{n-1}(\tilde f'\circ(f^\ell|W_{i_n\cdots i_0}))]      $$
The right-hand side is bounded by a sum of $n+1$ terms where ${\rm Var}$ is applied to $(\gamma-\pi_n\gamma)|W_{i_n\cdots i_0}$ or a factor $\tilde f'\circ(f^\ell|W_{i_n\cdots i_0}))$, and the other factors are bounded by their $||\cdot||_0$-norm.  Thus, using the hyperbolicity condition of Section 4, we have
$$	||({\cal L}_*^n-F_n)\gamma||      $$
$$	\le{\rm Var}(\gamma-\pi_n\gamma).A\alpha^n
+\sum_{\ell=0}^{n-1}\sum_{i_0\cdots i_n}||(\gamma-\pi_n\gamma)|W_{i_n\cdots i_0}||_0.
	A\alpha^\ell.{\rm Var}(\tilde f'|W_{i_{n-\ell}\cdots i_0}).A\alpha^{n-\ell-1}      $$
$$	\le2A\alpha^n{\rm Var}\,\gamma+nA^2\alpha^{n-1}{\rm Var}\,\tilde f'
	\sum_{i_0\cdots i_n}||(\gamma-\pi_n\gamma)|W_{i_n\cdots i_0}||_0      $$
$$	\le(2A+nA^2\alpha^{-1}{\rm Var}\,\tilde f')\alpha^n{\rm Var}\,\gamma
	\le(2A+nA^2\alpha^{-1}{\rm Var}\,\tilde f')\alpha^n||\gamma||      $$
so that
$$	||{\cal L}_*^n-F_n||\le(2A+nA^2\alpha^{-1}{\rm Var}\,\tilde f')\alpha^n      $$
and therefore ${\cal L}_*$ has essential spectral radius $\le\alpha<1$ on ${\cal H}$.  Suppose that there existed an eigenfunction $\gamma\in{\cal H}$ to the eigenvalue $1$ of ${\cal L}_*$; the fact that $\gamma$ is continuous and $\ne0$ on some $W_j$ would imply
$$	\int({\cal L}_*^n|\gamma|)(x)\,dx<\int|\gamma|(x)\,dx      $$
[because, for some $n$, ${\cal L}_*^n$ sends "mass" into $V_0$].  But this is in contradiction with
$$	\int|\gamma|(x)\,dx=\int|{\cal L}_*^n\gamma|(x)\,dx
	\le\int({\cal L}_*^n|\gamma|)(x)\,dx      $$
Therefore, $1$ cannot be an eigenvalue of ${\cal L}_*$, and there is $\gamma=(1-{\cal L}_*)^{-1}\eta\in H$ such that
$$	\gamma={\cal L}_*\gamma+\eta      $$
Since $\gamma^0$ satisfies the same equation in $L^1$, we have $\gamma^0-\gamma={\cal L}_*(\gamma^0-\gamma)$ hence $\gamma^0-\gamma=0$ by the same argument as above [$|\gamma^0-\gamma|$ is in $L^1$, with "mass" in some $V_\alpha$ because $H(u_1)$ has measure 0, and this is sent to $V_0$ by ${\cal L}_*^n$ for some $n$].  Thus $\gamma^0$ is continuous of bounded variation on the intervals $W_j$ for $j=1,\ldots,N$, and $\phi^0$ has bounded variation on $[a,b]$, with possible discontinuities only at $f^nu_1$ for $n=0,\ldots,N$, and $\phi^0(a)=\phi^0(b)=0$.  We have
$$	{\cal L}_{(1)}\phi^0-c_0^0\psi_0+\sum_{n=0}^\infty c_n^0\chi_n=\phi^0      $$
Therefore, hyperbolicity along the periodic orbit of $u_1$ shows that $\phi^0$ cannot have discontinuities, and this proves part (b) of Proposition 15.\qed
\medskip
	This also concludes the proof of Theorem 9.\qed
\bigskip\noindent
{\bf 16  Remarks.}
\medskip
	(a)  Theorem 9 shows that the density $\rho(x)$ of the unique a.c.i.m. $\rho(x)\,dx$ for $f$ can be written as the sum of spikes $\approx|x-f^nb|^{-1/2}\theta_n(x)$ (where $\theta_n$ vanishes unless $x>f^nb$ or $x<f^nb$) and a continuous background $\phi(x)$.  In fact, one can also write $\rho(x)$ as the sum of singular terms $\approx|x-f^nb|^{-1/2}\theta_n(x),|x-f^nb|^{1/2}\theta_n(x)$ and a background $\phi(x)$ which is now differentiable.  This result is discussed in Appendix A.  It seems clear that one could write $\rho(x)$ as a sum of terms $|x-f^nb|^{k/2}\theta_n(x)$ with $k=-1,1,\ldots,{2\ell-1\over2}$ and a background $\phi(x)$ of class $C^\ell$, but we have not written a proof of this.
\medskip
	(b)  Let $u\in(-\infty,u_1)\cup(u_1,v_1)\cup(v_2,u_2)\cup(u_2,\infty)$ and choose $w\in\{u_1,u_2,v_1,v_2\}$ such that $w$ is an endpoint of the interval containing $u$.  If $\pm(w-u)>0$ and $\theta_\pm$ is the characteristic function of $\{x:(w-x)(x-u)>0\}$ we define
$$	\psi_{(u\pm)}(x)={w-x\over w-u}\cdot|x-u|^{-1/2}\theta_\pm(x)      $$
or a similar expression with the same singularity at $u$, greater smoothness at $w$, and/or $\int\psi_{(u\pm)}=0$.  [Note that the $\psi_n$ are of this form].  Claim: if $u\in\tilde H$, there exists a unique $(\phi_\alpha)\in{\cal A}_1$ such that $\phi_\alpha=\psi_{(u\pm)}|V_\alpha$ for all $\alpha$; furthermore $||(\phi_\alpha)||_1$ has a bound independent of $u\pm$.  These results are proved in Appendix B (assuming $\gamma<\alpha^{-1/2}$).
\medskip
	Note that if $((\phi_\alpha),(c_n))\in{\cal A}$ and $c_0=c_1=0$, there is $(\tilde\phi_\alpha)\in{\cal A}_1$ such that $\Delta(\tilde\phi_\alpha)=w((\phi_\alpha),(c_n))$.  It seems thus that we might have replaced ${\cal A}$ by ${\cal A}_1$ in our earlier discussions.  However, separating the spikes $(c_n)$ from the background $(\phi_\alpha)$ was needed in the spectral study of ${\cal L}$.
\medskip
	(c)  The eigenvector $\Phi^0$ of ${\cal L}$ corresponding to the eigenvalue $1$ (with $w\Phi^0\ge0,\int w\Phi^0=1$) depends continuously on $f$.  To make sense of this statement we may consider a one-parameter family $(f_\kappa)$ such that $f_0=f$.  We let $H_\kappa,\tilde H_\kappa$ (hyperbolic sets) and ${\cal A}_{1\kappa}$ (Banach space) reduce to $H,\tilde H$ and ${\cal A}_1$ when $\kappa=0$.  We restrict $\kappa$ to a compact set $K$ such that $f_\kappa^3c_\kappa\in\tilde H_\kappa$ (where $c_\kappa$ is the critical point of $f_\kappa$).  The intervals $V_{\kappa\alpha}$ associated with $H_\kappa$ can be mapped to the $V_\alpha$ associated with $H$, providing an identification $\eta_\kappa:{\cal A}_{\kappa1}\to{\cal A}_1$.  There are natural definitions of ${\cal L}_\kappa:{\cal A}_{\kappa1}\oplus{\cal A}_2\to{\cal A}_{\kappa1}\oplus{\cal A}_2$ and the eigenvector $\Phi_\kappa^0$ reducing to ${\cal L}$ and $\Phi^0$ when $\kappa=0$.  We claim that $\kappa\mapsto\Phi_\kappa^\times=(\eta_\kappa,{\bf1})\Phi_\kappa^0$ is a continuous function $K\to{\cal A}_1\oplus{\cal A}_2$.  This result is proved in Appendix C.  It implies that, if $A$ is smooth, $\kappa\to\langle\Phi_{f_\kappa}^0,A\rangle$ is continuous on $K$.  The weight of the $n$-th spike is $C_0\prod_{k=1}^n|f'_\kappa(f_\kappa^{k-1}b_\kappa)|^{-1/2}$ and its speed is
$$	{d\over d\kappa}f_\kappa^nb_\kappa
	=\prod_{k=1}^nf'_\kappa(f_\kappa^{k-1}b_\kappa){db_\kappa\over d\kappa}
	+\sum_{\ell=1}^n\prod_{k=\ell+1}^nf'_\kappa(f_\kappa^{k-1}b_\kappa)
	f_\kappa^*(f_\kappa^{\ell-1}b_\kappa)\quad{\rm with}\quad 
	f_\kappa^*={df_\kappa\over d\kappa}      $$
The weight may be roughly estimated as $\sim\alpha^{n/2}$ and the speed as $\sim\alpha^{-n}$ for some $\alpha\in(0,1)$, suggesting that $\kappa\to\langle\Phi_{f_\kappa}^0,A\rangle$ is ${1\over2}$-H\"older on $K$.
\bigskip\noindent
{\bf 17  Informal study of the differentiability of $f\mapsto\langle\Phi_f^0,A\rangle$.}
\medskip
	Writing $\Phi_f^0$ instead of $\Phi^0$ we want to study the change of $\langle\Phi_f^0,A\rangle=\int dx\,(w\Phi_f^0)(x)A(x)$ when $f$ is replaced by $\hat f$ close to $f$ (and the critical orbit $\hat f^k\hat c$ for $k\ge3$ is in the perturbed hyperbolic set $\hat{\tilde H}$).  Writing $g={\rm id}-\hat f(\hat c)+f(c)$, we see that $\hat f$ is conjugate to $g\circ\hat f\circ g^{-1}$, which has maximum $f(c)$ at $g(\hat c)$.  With proper choice of the inverse $f^{-1}$ we have $f^{-1}\circ(g\circ\hat f\circ g^{-1})=h$ close to ${\rm id}$, hence $g\circ\hat f\circ g^{-1}=f\circ h$ and $(h\circ g)\circ\hat f\circ(h\circ g)^{-1}=h\circ f$, {\it i.e.}, $\hat f$ is conjugate to $h\circ f$ and we may write
$$	\langle\Phi_{\hat f}^0,A\rangle=\langle\Phi_{h\circ f}^0,A\circ h\circ g\rangle      $$
The differentiability of $\hat f\mapsto A\circ h\circ g$ is trivial, and we concentrate on the study of $h\mapsto\langle\Phi_{h\circ f}^0,A\rangle$.  Writing $h={\rm id}+X$, where $X$ is analytic, we see that the change $\delta(w\Phi_f^0)$ when $f$ is replaced by $({\rm id}+X)\circ f$ is, to first order in $X$, formally
$$	(1-{\cal L})^{-1}{\cal D}(-X\Phi_f^0)      $$
where ${\cal D}$ denotes differentiation.  [The above formula is standard first order perturbation calculation, and we have omitted the $w$ map from our formula].
\medskip
	Writing $\Phi_f^0=((\phi_\alpha^0),(C_n))$, we can identify ${\cal D}(-X((\phi_\alpha^0),0))$ with an element $\Phi^\times$ of ${\cal A}$ (so that $w\Phi^\times={\cal D}(Xw((\phi_\alpha^0),0))$ and $\int dx\,w\Phi^\times(x)=0$, use Appendix A) which is easy to study, and we are left to analyze the singular part ${\cal D}(-X(0,(C_n)))$.  To study this singular part we shall write $(0,(C_n))=\sum_{n=0}^\infty C_n\psi_{(f^nb)}$, and use the equivalence $\sim$ modulo the elements of ${\cal A}$.  We extend the domain of definition of ${\cal L}$ so that ${\cal L}\psi_{(u)}\sim|f'(u)|^{-1/2}\psi_{(fu)}$, where we use the notation $\psi_{(u\pm)}$ of Section 16(b), but omit the $\pm$, and we assume that $\int\psi_{(u)}=0$.  We have thus
$$	{\cal D}(-X(0,(C_n)))\sim-\sum_{n=0}^\infty C_nX(f^nb){\cal D}\psi_{(f^nb)}
	 \sim\sum_{n=0}^\infty C_nX(f^nb){d\over du}\psi_{(u)}\Big|_{u=f^nb}    $$
$$	=\sum_{n=0}^\infty C_nX(f^nb)[\prod_{k=0}^{n-1}f'(f^kb)]^{-1}{d\over db}\psi_{(f^nb)}
	\sim\sum_{n=0}^\infty X(f^nb)[\prod_{k=0}^{n-1}f'(f^kb)]^{-1}
	{d\over db}{\cal L}^nC_0\psi_{(b)}      $$
We may thus write (introducing $(1-\lambda{\cal L})^{-1}$ instead of $(1-{\cal L})^{-1}$)
$$	(1-\lambda{\cal L})^{-1}{\cal D}(-X(0,(C_n)))
\sim\sum_{n=0}^\infty X(f^nb)[\prod_{k=0}^{n-1}f'(f^kb)]^{-1}\lambda^{-n}
	{d\over db}(1-\lambda{\cal L})^{-1}(\lambda{\cal L})^nC_0\psi_{(b)}      $$
$$	=\sum_{n=0}^\infty X(f^nb)[\prod_{k=0}^{n-1}f'(f^kb)]^{-1}\lambda^{-n}
	{d\over db}(1-\lambda{\cal L})^{-1}C_0\psi_{(b)}-Z      $$
where
$$	Z=\sum_{n=0}^\infty X(f^nb)[\prod_{k=0}^{n-1}f'(f^kb)]^{-1}\lambda^{-n}
	{d\over db}\sum_{\ell=0}^{n-1}(\lambda{\cal L})^\ell C_0\psi_{(b)}      $$
$$	\sim\sum_{n=0}^\infty X(f^nb)[\prod_{k=0}^{n-1}f'(f^kb)]^{-1}
	\sum_{\ell=0}^{n-1}\lambda^{-n+\ell}|\prod_{k=0}^{\ell-1}f'(f^kb)|^{-1/2}
	{d\over db}C_0\psi_{(f^\ell b)}      $$
$$	=\sum_{n=0}^\infty X(f^nb)\sum_{\ell=0}^{n-1}\lambda^{-n+\ell}
	[\prod_{k=\ell}^{n-1}f'(f^kb)]^{-1}|\prod_{k=0}^{\ell-1}f'(f^kb)|^{-1/2}
	{d\over du}C_0\psi_{(u)}\Big|_{u=f^\ell b}    $$
$$	\sim-{\cal D}\sum_{n=0}^\infty X(f^nb)\sum_{\ell=0}^{n-1}\lambda^{-n+\ell}
	[\prod_{k=\ell}^{n-1}f'(f^kb)]^{-1}C_\ell\psi_\ell      $$
$$	=-{\cal D}\sum_{r=1}^\infty\sum_{\ell=0}^\infty X(f^{\ell+r}b)\lambda^{-r}
	[\prod_{k=0}^{r-1}f'(f^{\ell+k}b)]^{-1}C_\ell\psi_\ell      $$
$$	=-{\cal D}\sum_{\ell=0}^\infty C_\ell\psi_\ell\sum_{r=1}^\infty\lambda^{-r}
	[\prod_{k=0}^{r-1}f'(f^{\ell+k}b)]^{-1}X(f^{\ell+r}b)      $$
We have thus an (informal) proof of the following result
\medskip
{\sl	For $\ell=0,1,\ldots$, define
$$	F_\ell(X)=\sum_{n=1}^\infty\lambda^{-n}[\prod_{k=0}^{n-1}f'(f^{k+\ell}b)]^{-1}X(f^{n+\ell}b)  $$
which are holomorphic functions of $\lambda$ when $|\lambda|>\alpha$.  Then the {\rm susceptibility} function 
$$	\Psi(\lambda)=\langle(1-\lambda{\cal L})^{-1}{\cal D}(-X\Phi_f^0),A\rangle      $$
has the form
$$	\Psi(\lambda)\sim(X(b)+F_0(X))
	{d\over db}\langle(1-\lambda{\cal L})^{-1}C_0\psi_{(b)},A\rangle
	-\sum_{\ell=0}^\infty F_\ell(X)C_\ell\langle\psi_\ell,{\cal D}A\rangle      $$
The derivative ${d\over db}\langle(1-\lambda{\cal L})^{-1}C_0\psi_{(b)},A\rangle$
exists as a distribution, but is in principle a divergent quantity for given $b$.  The corresponding term disappears however if $X(b)+F_0(X)=0$, and we are then left with a finite expression, meromorphic in $\lambda$ for $\alpha<|\lambda|<\min(\beta^{-1},\alpha^{-1/2})$ and holomorphic when $\alpha<|\lambda|\le1$.}
\medskip
	Note that in writing the equivalence $\sim$ we have omitted terms with the singularities of $(1-\lambda{\cal L})^{-1}$; this explains the meromorphic contributions for $|\lambda|>1$.  The condition $X(b)+F_0(X)=0$ for $\lambda=1$ is known as {\it horizontality} (see the discussion in Section 19 below).
\bigskip\noindent
{\bf 18 A modified susceptibility function $\Psi(X,\lambda)$.}
\medskip
	At this point we extend the definition of the operator ${\cal L}$ to ${\cal L}^\sim$ acting on a larger space.  Remember that ${\cal L}$ was obtained from the transfer operator ${\cal L}_{(1)}$ by separating the spikes $\psi_n$ from the background in order to obtain better spectral properties.  We now also introduce derivatives $\psi'_n$ of spikes, so that the transfer operator sends $\psi'_n$ to
$$	{f'(f^nb)\over|f'(f^nb)|^{1/2}}\psi'_{n+1}+\hbox{a term in $w({\cal A}_1+{\cal A}_2)$}  $$
The coefficients of $\psi'_n$ form an element of ${\cal A}_3=\{(Y_n):||(Y_n)||_3=\sup_n\delta^n|Y_n|<\infty\}$.  We define ${\cal L}^\sim$ on ${\cal A}_1+{\cal A}_2+{\cal A}_3$ so that
$$	{\cal L}^\sim=\pmatrix{{\cal L}_0+{\cal L}_1&{\cal L}_2&{\cal L}_5\cr
	{\cal L}_3&{\cal L}_4&{\cal L}_6\cr
	0&0&{\cal L}_7\cr}      $$
where we omit the explicit definition of ${\cal L}_5$, ${\cal L}_6$, and let
$$	{\cal L}_7\Big({Z_n\over\prod_{k=0}^{n-1}|f'(f^kb)|^{1/2}}\Big)
	=\Big({\tilde Z_n\over\prod_{k=0}^{n-1}|f'(f^kb)|^{1/2}}\Big)      $$
with $\tilde Z_0=0$, $\tilde Z_n=f'(f^{n-1}b)Z_{n-1}$ for $n>0$.  Since
$$	\pmatrix{0&0&{\cal L}_5\cr0&0&{\cal L}_6\cr0&0&{\cal L}_7\cr}{\cal L}=0      $$
we have
$$	{\cal L}^{\sim n}={\cal L}^n
	+\sum_{k=1}^n{\cal L}^{k-1}({\cal L}_5+{\cal L}_6){\cal L}_7^{n-k}+{\cal L}_7^n      $$
and formally
$$	({\bf1}-\lambda{\cal L}^\sim)^{-1}=({\bf1}_{12}-\lambda{\cal L})^{-1}
	+({\bf1}_3-\lambda{\cal L}_7)^{-1}
	+({\bf1}_{12}-\lambda{\cal L})^{-1}\lambda({\cal L}_5+{\cal L}_6)
	({\bf1}_3-\lambda{\cal L}_7)^{-1}      $$
where ${\bf1}_{12}$ and ${\bf1}_3$ denote the identity on ${\cal A}_1\oplus{\cal A}_2$ and ${\cal A}_3$ respectively.
\medskip
	For $\lambda$ close to $1$, $({\bf1}_3-\lambda{\cal L}_7)^{-1}$ and thus $({\bf1}-\lambda{\cal L}^\sim)^{-1}$ are not well defined.  But there is a natural definition of a left inverse ${\cal L}_{7L}^{-1}$ of ${\cal L}_7$ where
$$	{\cal L}_{7L}^{-1}\Big({Z_n\over\prod_{k=0}^{n-1}|f'(f^kb)|^{1/2}}\Big)
	=\Big({\tilde Z_n\over\prod_{k=0}^{n-1}|f'(f^kb)|^{1/2}}\Big)      $$
with $\tilde Z_n=f'(f^nb)^{-1}Z_{n+1}$ for $n\ge0$.  The spectral radius of ${\cal L}_{7L}^{-1}$ is thus $\le\alpha^{1/2}/\delta$.  This gives natural left inverses
$$	({\bf1}_3-\lambda{\cal L}_7)_L^{-1}
	=-\sum_{n=1}^\infty\lambda^{-n}{\cal L}_{7L}^{-n}      $$
for $|\lambda|>\alpha^{1/2}/\delta$, and
$$	({\bf1}-\lambda{\cal L}^\sim)_L^{-1}=({\bf1}_{12}-\lambda{\cal L})^{-1}
	+({\bf1}_3-\lambda{\cal L}_7)_L^{-1}
	+({\bf1}_{12}-\lambda{\cal L})^{-1}\lambda({\cal L}_5+{\cal L}_6)
	({\bf1}_3-\lambda{\cal L}_7)_L^{-1}      $$
when $|\lambda|>\alpha^{1/2}/\delta$ and $({\bf1}_{12}-\lambda{\cal L})^{-1}$ exists.  This gives a modified susceptibility function
$$	\Psi_L(\lambda)
	=\langle({\bf1}-\lambda{\cal L}^\sim)_L^{-1}{\cal D}(-X\Phi_f^0),A\rangle      $$
meromorphic in $\lambda$ for $\alpha<|\lambda|<\min(\beta^{-1},\alpha^{-1/2})$ and holomorphic for $\alpha<|\lambda|\le1$.
\medskip
	Note that the ${\cal A}_3$ part of ${\cal D}(-X\Phi_f^0)$ is      
$$	(Y_n)
=\Big({-X(f^nb)\over{1\over2}|f''(c)|^{1/2}\prod_{k=0}^{n-1}|f'(f^kb)|^{1/2}}\Big)_{n\ge0}   $$
where $\sup_n|X(f^nb)|<\infty$.  Therefore, for small $|\lambda|$,
$$	({\bf1}_3-\lambda{\cal L}_7)^{-1}(Y_n)
	=\Big({-\sum_{k=0}^n\lambda^k(\prod_{\ell=1}^kf'(f^{n-\ell}b))X(f^{n-k}b)
	\over{1\over2}|f''(c)|^{1/2}\prod_{k=0}^{n-1}|f'(f^kb)|^{1/2}}\Big)_{n\ge0}      $$
because the right-hand side is in ${\cal A}_3$.  Note that the right-hand side is also in ${\cal A}_3$ under the condition
$$	\sum_{n=0}^\infty\lambda^{-n}(\prod_{k=0}^{n-1}f'(f^kb))^{-1}X(f^nb)=0\eqno{(*)}    $$
because this condition implies
$$	-\sum_{k=0}^n\lambda^{-k}(\prod_{\ell=0}^{k-1}f'(f^\ell b))^{-1}X(f^kb)
	=\sum_{k=n+1}^\infty\lambda^{-k}(\prod_{\ell=0}^{k-1}f'(f^\ell b))^{-1}X(f^kb)      $$
hence, multiplying by $\lambda^n\prod_{\ell=0}^{n-1}f'(f^\ell b)$,
$$		-\sum_{k=0}^n\lambda^{n-k}(\prod_{\ell=k}^{n-1}f'(f^\ell b))X(f^kb)
	=\sum_{k=n+1}^\infty\lambda^{n-k}(\prod_{\ell=n}^{k-1}f'(f^\ell b))^{-1}X(f^kb)      $$
or
$$	-\sum_{k=0}^n\lambda^k(\prod_{\ell=1}^kf'(f^{n-\ell}b))X(f^{n-k}b)
	=\sum_{k=1}^\infty\lambda^{-k}(\prod_{\ell=0}^{k-1}f'(f^{n+\ell}b))^{-1}X(f^{n+k}b)     $$
for each $n$, provided $|\lambda|>\alpha$.  We have proved that:
\medskip
	{\sl Under the condition $(*)$, a resummation of the series defining
$$	\langle({\bf1}-\lambda{\cal L}^\sim)^{-1}{\cal D}(-X\Phi_f^0),A\rangle      $$
yields $\Psi_L(\lambda)$.}  
\medskip
	It is then natural to define a modified susceptibility function $\Psi(X,\lambda)$ by
$$	(X,\lambda)\mapsto\Psi(X,\lambda)
	=\Psi_L(\lambda)\quad{\rm on}\quad\{(X,\lambda):\hbox{$(*)$ holds}\}      $$
Note that the left-hand side of $(*)$ is equal to the quantity $X(b)+F_0(X)$ met in Section 17, and that $(*)$ with $\lambda=1$ reduces to the horizontality condition.
\medskip
	We conclude with a rigorous result agreeing in part with the informal study in Section 17, in part with a conjecture of Baladi [3], Baladi and Smania [5].
\bigskip\noindent
{\bf 19 Theorem} (differentiability along topological conjugacy classes).
\medskip
	{\sl Let $f_\kappa=h_\kappa\circ f$ where the $h_\kappa$ are real analytic, depend smoothly on $\kappa$, and $f_\kappa^3c=\xi_\kappa f^3c$ identically in $\kappa$.  [This last condition expresses that $f_\kappa$ belongs to a conjugacy class, and $\xi_k:H\to H_\kappa$ is the conjugacy defined in Appendix C].  Then, if $A$ is smooth, $\kappa\mapsto\langle\Phi_\kappa^0,A\rangle=\int dx\,(w\Phi_{f_k}^0)(x)A(x)$ is continuously differentiable.  Furthermore
$$ {d\over d\kappa}\langle\Phi_{f_k}^0,A\rangle\Big|_{\kappa=0}=\Psi(X,1)      $$
where $\Psi(X,\lambda)$ is defined in Section 18 with $X={d\over d\kappa}h_\kappa|_{\kappa=0}$, and $\Psi(X,\lambda)$ is holomorphic for $\alpha<|\lambda|\le1$, meromorphic for $\alpha<|\lambda|<\min(\beta^{-1},\alpha^{-1/2})$.}
\medskip
	[The value $\kappa=0$ plays no special role, and is chosen for notational simplicity in the formulation of the theorem].
\medskip
	Our notion of topological conjugacy class is a special case of that discussed in [1].
\medskip
	Note that $\xi_0={\rm id}$, and that $\xi_\kappa$ depends differentiably on $\kappa$.
Since $f_\kappa^3c=\xi_\kappa f^3c$ and $f_\kappa\xi_\kappa=\xi_\kappa f$ on $H$, we have $f_\kappa^nc=\xi_\kappa f^nc$ for $n\ge3$ and by differentiation (writing $\xi'={d\over d\kappa}\xi_\kappa|_{\kappa=0}$):
$$	\sum_{k=1}^n[\prod_{\ell=k}^{n-1}f'(f^\ell c)]X(f^kc)=\xi'(f^nc)      $$
or
$$	\sum_{k=1}^n[\prod_{\ell=1}^{k-1}f'(f^\ell c)]^{-1}X(f^kc)
	=[\prod_{\ell=1}^{n-1}f'(f^\ell c)]^{-1}\xi'(f^nc)      $$
and letting $n\to\infty$:
$$	\sum_{k=1}^\infty[\prod_{\ell=1}^{k-1}f'(f^\ell c)]^{-1}X(f^kc)=0\qquad{\rm or}\qquad
	\sum_{n=0}^\infty[\prod_{k=0}^{n-1}f'(f^kb)]^{-1}X(f^nb)=0      $$
This is the horizontality condition derived much more generally in [1].
\medskip
	The proof of the theorem will be based on Appendices A, B, C, and use particularly the notation of Appendix C.  We write $\Phi_{f_\kappa}^0=\Phi_\kappa^0$ and recall that the expression
$$	\langle\Phi_\kappa^0,A\rangle_\kappa=\int dx\,(w_\kappa\Phi_\kappa^0)(x)A(x)
	=\sum_\alpha\int_{V_{\kappa\alpha}}\phi_{\kappa\alpha}^0A(x)dx
	+\sum_nc_{\kappa n}^0\int\psi_{\kappa n}(x)A(x)dx      $$
depends explicitly on the intervals $V_{\kappa\alpha}$ and the points $f_\kappa^kc$ for $k\ge1$.  We shall first prove the existence of ${d\over d\kappa}\langle\Phi_\kappa^0,A\rangle_\kappa|_{\kappa=0}=\lim_{\kappa\to0}{1\over\kappa}\int(w_\kappa\Phi_\kappa^0-w\Phi^0)A$ and give an expression involving only the intervals $V_\alpha$ and the points $f^kc$ (corresponding to $\kappa=0$).  Then we shall transform the expression obtained to the form $\Psi(X,1)$.
\medskip	
	Since the map $\xi_\kappa:H\to H_\kappa$ depends smoothly on $\kappa$ (in particular $f'_\kappa(f_\kappa^kb_\kappa)=f'_\kappa(\xi_\kappa f^kb)$ is continuous uniformly in $k$), it is easily seen that the operator ${\cal L}_\kappa^\times$ defined in Appendix C now depends continuously and even differentiably on $\kappa$.
\medskip
	We may write
$$	\langle\Phi_\kappa^0,A\rangle_\kappa
	=\sum_\alpha\int_{V_{\kappa\alpha}}\phi_{\kappa\alpha}^0(x)A(x)\,dx
	+\sum_nc_{\kappa n}^0\int\psi_{\kappa n}(x)A(x)\,dx      $$
$$	=\langle((\phi_{\kappa\alpha}^0),
	(c_{\kappa n}^0)),((A|V_{\kappa\alpha}),A)\rangle_\kappa      $$
$$	=\langle\Phi_\kappa^0,((A|V_{\kappa\alpha}),0)\rangle_\kappa
	+\langle\Phi_\kappa^0,(0,(c_{\kappa n}^0))\rangle_\kappa      $$
For notational simplicity we study the derivative of this quantity at $\kappa=0$ but the proof will show that the derivative depends continuously on $\kappa$.  We have
$$	{1\over\kappa}\Big[\langle\Phi_{\kappa}^0,A\rangle_{\kappa}
	-\langle\Phi^0,A\rangle\Big]=I+II      $$
where
$$	II={1\over\kappa}\sum_n\int
	[c_{\kappa n}^0\psi_{\kappa n}(x)-c_n^0\psi_n(x)]A(x)\,dx      $$
$$	\to\sum_n\int[{dc_{\kappa n}^0\over d\kappa}\psi_n(x)
	+c_n^0{d\over  d\kappa}\psi_{\kappa n}(x)]A(x)\,dx\Big|_{\kappa=0}      $$
[${d\over d\kappa}\psi_{\kappa n}$ is a distribution with singular part ${d\over d\kappa}|x-f_\kappa^nb_\kappa|^{-1/2}$; integrating by part over $x$, and using $f_ \kappa^nb_\kappa=\xi_\kappa f^nb$ for $k\ge2$, we see that the right-hand side makes sense, and is the limit of the left-hand side when $\kappa\to0$].
\medskip
	We also have
$$	\langle\Phi_\kappa^0,((A|V_{\kappa\alpha}),0)\rangle_\kappa
	=\langle\Phi_\kappa^\times,((A_{\kappa\alpha}),0)\rangle      $$
where $A_{\kappa\alpha}=(A|V_{\kappa\alpha})\circ\tilde\eta_{\kappa\alpha}^{-1}$, so that
$$	I=\langle{\Phi_{\kappa}^\times-\Phi_0^\times\over\kappa},
	((A_{\kappa\alpha}),0)\rangle+\langle\Phi_0^\times,
	(({A_{\kappa\alpha}-A_{0\alpha}\over\kappa}),0)\rangle      $$
and the second term is readily seen to tend to a limit when $\kappa\to0$.  In the first term remember that for $\kappa=0$ we have $\Phi_\kappa^\times=\Phi_0^\times=\Phi^0$, and ${\cal L}_\kappa^\times={\cal L}_0^\times={\cal L}$.  Also
$$	({\bf 1}-{\cal L})(\Phi_{\kappa}^\times-\Phi_0^\times)
	=({\cal L}_{\kappa}^\times-{\cal L}_0^\times)\Phi_\kappa^\times      $$
hence
$$	\Phi_\kappa^\times-\Phi_0^\times=({\bf 1}-{\cal L})^{-1}
	({\cal L}_\kappa^\times-{\cal L}_0^\times)\Phi_\kappa^\times      $$
Since $({\bf 1}-{\cal L})^{-1}$ is bounded and $\kappa\mapsto{\cal L}_\kappa^\times$ differentiable, we have
$$	\langle{\Phi_\kappa^\times-\Phi_0^\times\over\kappa},
	((A_{\kappa\alpha}),0)\rangle\to\langle({\bf 1}-{\cal L})^{-1}({d\over d\kappa}
	{\cal L}_\kappa^\times\Big|_{\kappa=0})\Phi^0,((A_{0\alpha}),0)\rangle      $$
when $\kappa\to0$, proving that $\kappa\mapsto\langle\Phi_\kappa^0,A\rangle$ is differentiable.
\medskip
	If we replace in the above calculation the Banach space ${\cal A}={\cal A}_1\oplus{\cal A}_2$ by ${\cal A}'={\cal A}'_1\oplus{\cal A}'_2$ as in Appendix A, we obtain an expression of ${d\over d\kappa}\langle\Phi_\kappa^0,A\rangle_\kappa|_{\kappa=0}$ that can be re-expressed in terms of the ${\psi}'_n$, ${\psi}_n$ and an element of ${\cal A}_1$.  We may thus write
$$	{d\over d\kappa}\langle\Phi_\kappa^0,A\rangle_\kappa\Big|_{\kappa=0}
	=\langle\tilde\Phi,A\rangle^\sim      $$
where $\tilde\Phi\in{\cal A}_1\oplus{\cal A}_2\oplus{\cal A}_3$.  The part $\tilde\Phi_3$ of $\tilde\Phi$ in ${\cal A}_3$ is uniquely determined by $A\mapsto\langle\tilde\Phi,A\rangle^\sim$; the calculation of II above shows that $n$-th component of $\tilde\Phi_3$ is
$$	-{d\over d\kappa}f_\kappa^nb_\kappa\Big|_{\kappa=0}c_n^0
	=-{d\over d\kappa}f_\kappa^{n+1}c\Big|_{\kappa=0}c_n^0      $$
$$	=-\sum_{k=1}^{n+1}X(f^kc)\big(\prod_{\ell=k}^nf'(f^\ell c)\big)c_n^0
	=-\sum_{k=0}^nX(f^kb)\big(\prod_{\ell=k}^{n-1}f'(f^\ell b)\big)c_n^0      $$
and as a result
$$	({\bf1}-{\cal L}_7)\tilde\Phi_3=(-X(f^nb)C_n^0)_{n\ge0}      $$
$$	\tilde\Phi_3=({\bf1}-{\cal L}_7)_L^{-1}(-X(f^nb)C_n^0)_{n\ge0}      $$
The part $\Phi^*$ of $\tilde\Phi$ in ${\cal A}_1\oplus{\cal A}_2$ is not uniquely determined (because of the ambiguity discussed in Appendix B); this part satisfies $\int w\Phi^*=0$.
\medskip
	If ${\cal L}_{(1)\kappa}$ is the transfer operator corresponding to $f_\kappa$, we have ${\cal L}_{(1)\kappa}w_\kappa\Phi_\kappa^0=w_\kappa\Phi_\kappa^0$, hence
$$	({\bf1}-{\cal L}_{(1)})(w_\kappa\Phi_\kappa^0-w\Phi^0)
	=({\cal L}_{(1)\kappa}-{\cal L}_{(1)})w_\kappa\Phi_\kappa^0      $$
Therefore (using the fact that we may let ${\cal L}_{(1)}$ act on $A$) we have
$$	\langle({\bf1}-{\cal L}^\sim)\tilde\Phi,A\rangle^\sim
	=\lim_{\kappa\to0}\int A{1\over\kappa}({\bf1}-{\cal L}_{(1)})
	(w_\kappa\Phi_\kappa^0-w\Phi^0)      $$
$$	=\lim_{\kappa\to0}\int A{1\over\kappa}({\cal L}_{(1)\kappa}-{\cal L}_{(1)})
	w_\kappa\Phi_\kappa^0
	=\lim_{\kappa\to0}\int A{1\over\kappa}({\rm id}^*-h_{-\kappa}^*)
	w_\kappa\Phi_\kappa^0
	=\lim_{\kappa\to0}\int A{1\over\kappa}(h_\kappa^*-{\rm id}^*)w\Phi^0      $$
where $h^*$ denotes the direct image of a measure (here a L$^1$ function) under $h$, and the last equality uses the existence of a continuous derivative for $\kappa\mapsto\langle\Phi_\kappa^0,A\rangle_\kappa$.  According to Appendix A we may write $w\Phi^0$ as a sum of terms $C_n^{(0)}\psi_n^{(0)}$, $C_n^{(1)}\psi_n^{(1)}$, and a differentiable background.  Corresponding to this we may identify $\lim_{\kappa\to0}{1\over\kappa}(h_\kappa^*-{\rm id}^*)\Phi^0$ with a naturally defined element ${\cal D}(-X\Phi^0)$ of ${\cal A}_1\oplus{\cal A}_2\oplus{\cal A}_3$, where ${\cal D}$ denotes differentiation.  We write ${\cal D}(-X\Phi^0)=(D^*,D_3)$ with $D^*\in{\cal A}_1\oplus{\cal A}_2,D_3\in{\cal A}_3$.  Since the coefficient of $\psi'_n$ in ${\cal D}(-X\Phi^0)$ is $-X(f^nb)c_n^0$, we have $D_3=({\bf1}-{\cal L}_7)\tilde\Phi_3$. With $\tilde\Phi=(\Phi^*,\tilde\Phi_3)$ we have thus
$$	\langle({\bf1}-{\cal L}^\sim)(\Phi^*,\tilde\Phi_3),A\rangle^\sim
	=\langle{\cal D}(-X\Phi^0),A\rangle^\sim      $$
and
$$	\langle({\bf1}-{\cal L})\Phi^*,A\rangle
	=\langle{\cal D}(-X\Phi^0)-({\bf1}-{\cal L}^\sim)(0,\tilde\Phi_3),A\rangle      $$
In particular $\int w[{\cal D}(-X\Phi^0)-({\bf1}-{\cal L}^\sim)(0,\tilde\Phi_3)]=0$ and we may define
$$	\Phi=({\bf1}-{\cal L})^{-1}[{\cal D}(-X\Phi^0)-({\bf1}-{\cal L}^\sim)(0,\tilde\Phi_3)]
	\in{\cal A}      $$
We have then $\langle({\bf1}-{\cal L})(\Phi^*-\Phi),A\rangle=0$, hence
$$	w({\bf1}-{\cal L})(\Phi^*-\Phi)=0      $$
hence
$$	w(\Phi^*-\Phi)={\cal L}_{(1)}w(\Phi^*-\Phi)      $$
with $\int w(\Phi^*-\Phi)=0$, so that $w(\Phi^*-\Phi)=0$, and
$$	\langle\Phi^*,A\rangle=\langle\Phi,A\rangle=\langle({\bf1}-{\cal L})^{-1}
	[{\cal D}(-X\Phi^0)-({\bf1}-{\cal L}^\sim)(0,\tilde\Phi_3)],A\rangle      $$
$$	=\langle({\bf1}-{\cal L})^{-1}[D^*+{\cal L}_5+{\cal L}_6)\tilde\Phi_3],A\rangle
	=\langle({\bf1}-{\cal L})^{-1}
	[D^*+{\cal L}_5+{\cal L}_6)({\bf1}-{\cal L}_7)_L^{-1}D_3],A\rangle      $$
$$	=\langle({\bf1}-{\cal L}^\sim)_L^{-1}(D^*,D_3),A\rangle
	-\langle({\bf1}-{\cal L}_7)_L^{-1}D_3,A\rangle
	=\Psi(X,1)-\langle(0,\tilde\Phi_3),A\rangle^\sim      $$
so that finally
$$	{d\over d\kappa}\langle\Phi_\kappa^0,A\rangle_\kappa|_{\kappa=0}
	=\langle\tilde\Phi,A\rangle^\sim=\Psi(X,1)      $$
as announced.\qed
\medskip
	Note that in [5], Baladi and Smania study (in the case of piecewise expanding maps) the more difficult problem of differentiability in horizontal directions ({\it i.e.}, directions tangent to a topological class).  It appears likely that this could be done here also (as conjectured in [5]) , but we have not tried to do so.
\bigskip\noindent
{\bf 20  Discussion.}
\medskip
	The codimension 1 condition $X(b)+F_0(X)=0$ for $\lambda=1$ expresses that $X$ is a {\it horizontal} perturbation, which means that it is tangent to a topological class of unimodal maps (see [1] and references given there).  In our case, a family $(f_\kappa)$ is in a topological conjugacy class if $f_\kappa^3c_\kappa=\xi_\kappa f^3c$ in the notation of Appendix C.  The informal result obtained in Section 17 and the formal proof of differentiability along a topological conjugacy class given by Theorem 19 support the conjecture by Baladi and Smania [5] that the map $f\mapsto\langle\Phi_f^0,A\rangle$ is differentiable (in the sense of Whitney) in horizontal directions, {\it i.e.}, along a curve tangent to a topological conjugacy class.  Our theorem 19 also relates the derivative along a topological conjugacy class to a naturally defined susceptibility function.  It seems unlikely that a derivative (in the sense of Whitney) exists in nonhorizontal directions.  Note however that if $f\mapsto\langle\Phi_f^0,A\rangle$ is nondifferentiable, it will be in a mild way: the "nondifferentiable" contribution to $\Psi(\lambda)$ is, as we saw above, proportional to
$$	{d\over db}\langle(1-\lambda{\cal L})^{-1}\psi_{(b)},A\rangle
	\sim\sum_n\lambda^n{d\over db}\langle{\cal L}^n\psi_{(b)},A\rangle      $$
where $\langle{\cal L}^n\psi_{(b)},A\rangle$ decreases exponentially with $n$, while ${d\over db}\langle{\cal L}^n\psi_{(b)},A\rangle$ increases exponentially.  Therefore, if one does not see the small scale fluctuations of $b\mapsto\langle(1-\lambda{\cal L})^{-1}\psi_{(b)},A\rangle$, this function will seem differentiable.  But  the singularities with respect to $\lambda$ (with $|\lambda|<1$) may remain visible.  In conclusion, a physicist may see singularities with respect to $\lambda$ of a derivative (with respect to $f$ or $b$) while this derivative may not exist for a mathematician.
\vfill\eject\noindent
{\bf A  Appendix} (proof of Remark 16(a)).
\medskip
	We return to the analysis in Section 10, and note that by an analytic change of variable $x\mapsto\xi(x)$ we can get $y=fx=b-\xi^2$ [we have indeed $b-y=A(x-c)^2(1+\beta(x).(x-c))$ with $\beta$ analytic, and we can take $\xi=(x-c)A^{1/2}(1+\beta(x).(x-c))^{1\over2}$].  Write $\rho(x)\,dx=\tilde\rho(\xi)\,d\xi$ (where $\tilde\rho$ is analytic near $0$).  The density of the image $\delta(y)\,dy$ by $f$ of $\rho(x)\,dx=\tilde\rho(\xi)\,d\xi$ is, near $b$,
$$	\delta(y)={1\over2\sqrt{y-b}}(\tilde\rho(\sqrt{y-b})+\tilde\rho(-\sqrt{y-b}))
	={\hat\rho(y-b)\over\sqrt{y-b}}      $$
where $\hat\rho$ is analytic near $0$.  Therefore, near $b$,
$$	\delta(x)={U\over\sqrt{b-x}}+U'\sqrt{b-x}+\ldots      $$
where $U=\rho(c)/\sqrt{A}$, and $U'$ is linear in $\rho(c),\rho'(c),\rho''(c)$ with coefficients depending on the derivatives of $f$ at $c$.  Near $a$ we find
$$	\delta(x)=U|f'(b)|^{-1/2}{1\over\sqrt{x-a}}
	+(U'|f'(b)|^{-3/2}-{3\over4}Uf''(b)|f'(b)|^{-5/2})\sqrt{x-a}      $$
Writing $s_n=-{\rm sgn}\prod_{k=0}^{n-1}f'(f^kb)$, $t_n=|\prod_{k=0}^{n-1}f'(f^kb)|^{-1/2}$, we claim that near $f^nb$ we have a singularity given for $s_n(x-f^nb)<0$ by $0$, and for $s_n(x-f^nb)>0$ by 
$$	\delta(x)={Ut_n\over\sqrt{s_n(x-f^nb)}}+(U't_n^3-{3\over4}Ut_n
\sum_{k=0}^{n-1}s_{k+1}{f''(f^kb)\over|f'(f^kb)|}{t_n^2\over t_k^2})\sqrt{s_n(x-f^nb)}      $$
[to prove this we use induction on $n$, and the fact that, when $f:x\mapsto y$ for $x$ close to $f^nb$ we have:
$$	s_n(x-f^nb)={s_{n+1}(y-f^{n+1}b)\over|f'(f^nb)|}
	[1-{f''(f^nb)\over2|f'(f^nb)|^2}(y-f^{n+1}b)]      $$
$$dx={dy\over|f'(f^nb)|}[1-{f''(f^nb)\over|f'(f^nb)|^2}(y-f^{n+1}b)]\qquad\qquad].$$
Define now
$$	\psi_n^{(0)}(x)
	=\big(1-({x-f^nb\over w_n-f^nb})^2\big){\theta_n(x)\over\sqrt{s_n(x-f^nb)}}       $$
$$	\psi_n^{(1)}(x)
	=\big(1-({x-f^nb\over w_n-f^nb})^2\big)\theta_n(x)\sqrt{s_n(x-f^nb)}       $$
for $s_n(x-f^nb)>0$, $0$ otherwise.  Then, the expected singularity of $\delta$ near $f^nb$ is given by
$$	Ut_n\psi_n^{(0)}+(U't_n^3-{3\over4}Ut_n\sum_{k=0}^{n-1}s_{k+1}{f''(f^kb)\over
|f'(f^kb)|}	{t_n^2\over t_k^2})\psi_n^{(1)}=C_n^{(0)}\psi_n^{(0)}+C_n^{(1)}\psi_n^{(1)}   $$
where $C_0^{(0)}=U$, $C_0^{(1)}=U'$, and
$$	C_{n+1}^{(0)}=|f'(f^nb)|^{-1/2}C_n^{(0)}      $$
$$	C_{n+1}^{(1)}=|f'(f^nb)|^{-3/2}C_n^{(1)}-{3\over4}s_{n+1}
	|f'(f^nb)|^{-5/2}f''(f^nb)C_n^{(0)}  $$
$$	=|f'(f^nb)|^{-3/2}(C_n^{(1)}-{3\over4}s_{n+1}{f''(f^nb)\over|f'(f^nb)|}C_n^{(0)})  $$
Let
$$	f(\psi_n^{(0)}(x)\,dx)=\tilde\psi_{n+1}^{(0)}(x)\,dx\qquad,
	\qquad f(\psi_n^{(1)}(x)\,dx)=\tilde\psi_{n+1}^{(1)}(x)\,dx      $$
and write	
$$	\tilde\psi_{n+1}^{(0)}=|f'(f^nb)|^{-1/2}\psi_{n+1}^{(0)}
	-{3\over4}s_{n+1}|f'(f^nb)|^{-5/2}f''(f^nb)\psi_{n+1}^{(1)}+\chi_n^{(0)}      $$
$$	\tilde\psi_{n+1}^{(1)}=|f'(f^nb)|^{-3/2}\psi_{n+1}^{(1)}+\chi_n^{(1)}      $$
The density of $f(C_n^{(0)}\psi_n^{(0)}(x)\,dx+C_n^{(1)}\psi_n^{(1)}(x)\,dx)$ is then
$$	C_{n+1}^{(0)}\psi_{n+1}^{(0)}+C_{n+1}^{(1)}\psi_{n+1}^{(1)}
	+C_n^{(0)}\chi_n^{(0)}+C_n^{(1)}\chi_n^{(1)}      $$
The functions $\chi_n^{(0)},\chi_n^{(1)}$ have been constructed such that they and their first derivatives $\chi_n^{(0)\prime},\chi_n^{(1)\prime}$ have the properties of Lemma 11.  Namely, $\chi_n^{(0)},\chi_n^{(1)},\chi_n^{(0)\prime},\chi_n^{(1)\prime}$ are continuous with bounded variation on $[a,b]$ uniformly in $n$, they vanish at $a,b$, and if $n\ge1$ they extend to holomorphic functions on the appropriate $D_\alpha$, with uniform bounds.
\medskip
	Let ${\cal A}_1'\subset {\cal A}_1$ consist of the $(\phi_\alpha)$ such that the derivatives $\phi_{-1}',\phi_{-2}'$ of $\phi_{-1},\phi_{-2}$ vanish at $\pi_b^{-1}b$ and $\pi_a^{-1}a$ respectively.  Let also ${\cal A}_2'$ consist of the sequences $(c_n^{(0)},c_n^{(1)})$, with $c_n^{(0)},c_n^{(1)}\in{\bf C}$, $n=0,1,\ldots$ such that
$$	||(c_n^{(0)},c_n^{(1)})||_2'=\sup_{n\ge0}\delta^n(|c_n^{(0)}|+|c_n^{(1)}|)<\infty      $$
If $\Phi'=((\phi_\alpha),(c_n^{(0)},c_n^{(1)}))\in{\cal A}'={\cal A}_1'\oplus{\cal A}_2'$ we let $||\Phi'||'=||(\phi_\alpha)||_1+||(c_n^{(0)},c_n^{(1)})||_2'$, making ${\cal A}'$ into a Banach space.  We may now proceed as in Section 12, replacing ${\cal A}$ by ${\cal A}'$, and defining ${\cal L}':{\cal A}'\mapsto{\cal A}'$ in a way similar to ${\cal L}:{\cal A}\mapsto{\cal A}$, but with (ii), (v), (vi) replaced as follows:
\medskip
	(ii) $\displaystyle\phi_0\Rightarrow\Big((\hat c_0^{(0)},\hat c_0^{(1)})=(U,U')\,,\,\hat\phi_{-1}=\pm{\phi_0\over|f'|}\circ\tilde f_{-1}^{-1}-U(\pm{1\over2}\psi_0^{(0)}\circ\pi_b)-U'(\pm{1\over2}\psi_0^{(1)}\circ\pi_b)\Big)$
so that $\hat\phi_{-1}$ is holomorphic in $\pi_b^{-1}D_{-1}$ with vanishing derivative at $\pi_b^{-1}b$
\medskip
   (v) $\displaystyle(c_0^{(0)},c_0^{(1)})\Rightarrow\Big((\hat c_1^{(0)},\hat c_1^{(1)})=(|f'(b)|^{-1/2}c_0^{(0)},|f'(b)|^{-3/2}c_0^{(1)}-{3\over4}|f'(b)|^{-5/2}f''(b)c_0^{(0)})$, $\displaystyle\chi_0=\pm{1\over2}c_0^{(0)}({\psi_0^{(0)}\over|f'|}\circ\pi_b\circ\tilde f_{-2}^{-1}-|f'(b)|^{-1/2}\psi_1^{(0)}\circ\pi_a+{3\over4}|f'(b)|^{-5/2}f''(b)\psi_1^{(1)}\circ\pi_a)\qquad\qquad$\break $\displaystyle\pm{1\over2}c_0^{(1)}({\psi_0^{(1)}\over|f'|}\circ\pi_b\circ\tilde f_{-2}^{-1}-|f'(b)|^{-3/2}\psi_1^{(1)}\circ\pi_a)\hbox{ in $\pi_a^{-1}D_{-2}$}\Big)$
\medskip          
	(vi) $\displaystyle(c_n^{(0)},c_n^{(1)})\Rightarrow\qquad\qquad\qquad\qquad\qquad\qquad\qquad\qquad\qquad\qquad\qquad\qquad\qquad\qquad\qquad\qquad\qquad\qquad\qquad\qquad$\break $\Big((\hat c_{n+1}^{(0)},\hat c_{n+1}^{(1)})=(|f'(f^nb)|^{-1/2}c_n^{(0)},|f'(f^nb)|^{-3/2}c_n^{(1)}-{3\over4}s_{n+1}|f'(f^nb)|^{-5/2}f''(f^nb)c_n^{(0)})$,$\qquad\qquad\qquad\qquad\qquad\qquad\qquad\qquad\qquad\qquad\qquad\qquad\qquad$ \break$\displaystyle\chi_{n\alpha}=c_n^{(0)}[{\psi_n^{(0)}\over|f'|}\circ f_n^{-1}-|f'(f^nb)|^{-1/2}\psi_{n+1}^{(0)}+{3\over4}s_{n+1}|f'(f^nb)|^{-5/2}f''(f^nb)|\psi_{n+1}^{(1)}]\qquad\qquad\qquad\qquad\qquad\qquad\qquad\qquad\qquad\qquad\qquad\qquad\qquad$\break$\displaystyle+c_n^{(1)}[{\psi_n^{(1)}\over|f'|}\circ f_n^{-1}-|f'(f^nb)|^{-3/2}\psi_{n+1}^{(1)}]\qquad\hbox{in $D_\alpha$ if $V_\alpha\subset\{x:\theta_n(f_n^{-1}x)>0\}$, $0$ otherwise}\Big)$\break if $n\ge1$.
\medskip\noindent 
We write then
$$	{\cal L}'\Phi'=\tilde\Phi'=((\tilde\phi_\alpha),(\tilde c_n^{(0)},\tilde c_n^{(1)}))      $$
where
$$	\tilde\phi_{-2}=\hat\phi_{-2}+\chi_0\qquad,\qquad\tilde\phi_{-1}=\hat\phi_{-1}      $$
$$	\tilde\phi_\alpha=\sum_{\beta:fV_\beta=V_\alpha}\hat\phi_{\beta\alpha}
	+\hat\phi_\alpha+\sum_{n\ge1}\chi_{n\alpha}\qquad\hbox{if order $\alpha\ge0$}   $$
$$	(\tilde c_n^{(0)},\tilde c_n^{(1)})
	=(\hat c_n^{(0)},\hat c_n^{(1)})\qquad\hbox{for $n\ge0$}      $$
\indent
	With the above definitions and assumptions we find, by analogy with Theorem 13, that ${\cal L}':{\cal A}'\to{\cal A}'$ has essential spectral radius $\le\max(\gamma^{-1},\delta\alpha^{1/2})$.  There is (see Proposition 15) a simple eigenvalue $1$, and the rest of the spectrum has radius $<1$.  It is convenient to denote by $\Phi^0=((\phi_\alpha^0),(c_n^{0(0)},c_n^{0(1)}))$ the eigenfunction to the eigenvalue $1$.  We find again that $\phi^0=\Delta(\phi_\alpha^0)$ is continuous, of bounded variation, and satisfies $\phi^0(a)=\phi^0(b)=0$, but we can say more.  Using the notation in the proof of Proposition 15, we have again
$$	\gamma_j^0=\sum_k{\cal L}_{jk}\gamma_k^0+\eta_j      $$
with $\eta_j=\sum_{n=0}^\infty\eta_{jn}$, but now $\eta_{jn}=c_n^{0(0)}\chi_n^{(0)}+c_n^{0(1)}\chi_n^{(1)}|W_j$ for $n\ge1$, so that the $\eta_j$ have derivatives $\eta'_j\in{\cal H}_j$.  The derivatives $\gamma_j^{0\prime}$of the $\gamma_j^0$ are measures satisfying
$$	\gamma_j^{0\prime}=\sum_k{\cal L}_{jk}^\prime\gamma_k^{0\prime}+\eta_j^*     $$
The operator ${\cal L}_{jk}^\prime$ has the same form as ${\cal L}_{jk}$, but with an extra denominator $f'\circ(f^{-1})_{kj}$, and therefore ${\cal L}_*^\prime=({\cal L}_{jk}^\prime)$ acting on measures has spectral radius $\le\alpha<1$.  The term $\eta_j^*$ is the sum of $\eta_j^\prime$ and a term$\sum_k{\cal L}_{kj}^\sim\gamma_k^0$ where ${\cal L}_{kj}^\sim$ involves the derivative of $|f'\circ(f^{-1})_{kj}|^{-1}$ so that $\eta_j^*\in{\cal H}_j$.  The operator ${\cal L}_*^\prime$ also maps ${\cal H}$ to ${\cal H}$ and, by the same argument as for ${\cal L}_*$, has essential spectral radius $<1$ on ${\cal H}$.  Furthermore, $1$ cannot be an eigenvalue since ${\cal L}_*^\prime$ has spectral radius $<1$ on measures.  It follows that $(\gamma^{0\prime})=(\gamma_j^{0\prime})=(1-{\cal L}_*^\prime)^{-1}(\eta_j^*)\in{\cal H}$.  Therefore, the derivative $\phi^{0\prime}$ of $\phi^0$ may have discontinuities only on the orbit of $u_1$, and hyperbolicity again shows that this cannot happen.  In conclusion, $\phi^0$ and its derivative $\phi^{0\prime}$ are both of bounded variation, continuous, and vanishing at $a,b$.
\medskip
	A discussion similar to the above shows that the equation $\gamma=(1-{\cal L}_*^\prime)^{-1}\eta^*$ also defines $\gamma$ with finite norm in ${\cal A}_1$, and this $\gamma$ must coincide with $(\gamma^{0\prime})$ as a measure.  Therefore the family of derivatives $(\phi_\alpha^{0\prime})$ is an element of ${\cal A}_1$.  [For simplicity, we have written   
$\phi_{-1}^{0\prime}$, $\phi_{-2}^{0\prime}$ for the functions which, under application of $\Delta$, give the derivative of $\Delta\phi_{-1}^0$, $\Delta\phi_{-2}^0$].\qed
\vfill\eject\noindent
{\bf B  Appendix} (proof of Remark 16(b)).
\medskip
	If $u\in\tilde H$ and $\psi_{(u\pm)}$ is defined as in Remark 16(b), we want to show that there is a unique $(\phi_\alpha)$ in ${\cal A}_1$ such that $\phi_\alpha=\psi_{(u\pm)}|V_\alpha$ for all $\alpha$.  Furthermore $||(\phi_\alpha)||_1$ is bounded uniformly for $u\in\tilde H$, provided we assume $1<\gamma<\min(\beta^{-1},\alpha^{-1/2})$.
\medskip
	Note that uniqueness is automatic, and that $\phi_\alpha=0$ unless order $V_\alpha>0$.  Omitting the $\pm$ we let
$$	f(\psi_{(f^nu)}(x)\,dx)=[|f'(f^nu)|^{-1/2}\psi_{(f^{n+1}u)}(x)+\chi_{(f^nu)}(x)]\,dx       $$
For $n\ge0$ there is a unique $\omega_{un}$ such that $f^{n+1}(\omega_{un}\,dx)=\prod_{k=0}^{n-1}|f'(f^ku)|^{-1/2}\chi_{(f^nu)}\,dx$ and $[f^ku-c]\times[$supp $f^k(\omega_{un}(x)\,dx)-c]>0$ for $0\le k\le n$.  Furthermore $\psi_{(u)}=\sum_{n=0}^\infty\omega_{un}$ where the sum restricted to each $V_\alpha$ is finite.  If $[\chi_{(f^nu)}]$ denotes the element of ${\cal A}_1$ corresponding to $\chi_{(f^nu)}$, we find that $||[\chi_{(f^nu)}]||_1$ is bounded uniformly in $n$ and $u$.  Also note that we obtain $\omega_{un}$ from $\prod_{k=0}^{n-1}|f'(f^ku)|^{-1/2}\chi_{(f^nu)}$ by multiplying with $\prod_{k=0}^{n-1}|f'(f^ku)|$ (up to a factor bounded uniformly in $n$ because of hyperbolicity) and composing with $f^{n+1}$ (restricted to a small interval $J$ such that $f^{n+1}|J$ is invertible).  We have thus
$$	||[\omega_{un}]||_1\le{\rm const}\,\gamma^n\prod_{k=0}^{n-1}|f'(f^ku)|^{-1/2}      $$
where $[\omega_{un}]$ is the element of ${\cal A}_1$ corresponding to $\omega_{un}$ [This is because the replacement of $|V_\alpha|$ by $|(f|J)^{-n-1}V_\alpha|$ in the definition of $||\cdot||_1$ is compensated up to a multiplicative constant by the factor $\prod_{k=0}^{n-1}|f'(f^ku)|$ ].  Thus
$$	||[\omega_{un}]||_1\le{\rm const}\,(\gamma\alpha^{1/2})^n      $$
Since $\gamma<\alpha^{-1/2}$ we find that $\sum_n||[\omega_{un}]||_1<$ constant independent of $u$.  Therefore, since $(\phi_\alpha)=\sum_n[\omega_{un}]$, we see that $||(\phi_\alpha)||_1$ is bounded independently of $u$.\qed
\vfill\eject\noindent
{\bf C  Appendix} (proof of Remark 16(c)).
\medskip
	We consider a one-parameter family $(f_\kappa)$ of maps, reducing to $f=f_0$ for $\kappa=0$.  We assume that $(\kappa,x)\mapsto f_\kappa x$ is real-analytic.  For $\kappa$ close to $0$, $f_\kappa$ has a critical point $c_\kappa$ and maps $[a_\kappa,b_\kappa]$ to itself, with $b_\kappa=f_\kappa c_\kappa,a_\kappa=f_\kappa^2 c_\kappa$.  There is (by hyperbolicity of $H$ with respect to $f$) a homeomorphism $\xi_\kappa:H\to H_\kappa$ where $H_\kappa$ is an $f_\kappa$-invariant hyperbolic set for $f_\kappa$ and $f_\kappa\circ\xi_\kappa=\xi_\kappa\circ f$ on $H$.  We shall consider a compact set $K$ of values of $\kappa$ such that $f_\kappa a_\kappa\in\tilde H_\kappa$; we let $K\ni0$, $K$ of small diameter, and assume now $\kappa\in K$.  We may in a natural way define a Banach space ${\cal A}_\kappa={\cal A}_{\kappa1}\oplus{\cal A}_2$ and an operator ${\cal L}_\kappa:{\cal A}_\kappa\to{\cal A}_\kappa$ associated with $f_\kappa$ so that ${\cal A}_\kappa,{\cal L}_\kappa$ reduce to ${\cal A},{\cal L}$ for $\kappa=0$. Note that, since $\kappa\in K$ is close to $0$, we may assume that the constants $A,\alpha$ in the definition (Section 4) of hyperbolicity, and the constants $B,\beta$ (Section 7) are uniform in $\kappa$.
\medskip
	Let $\eta_{\kappa,-2}$ be a biholomorphic map of the complex neighborhood $D_{-2}$ of $[a,u_1]$ to the complex neighborhood $D_{\kappa,-2}$ of the corresponding interval $[a_\kappa,u_{\kappa1}]$, and lift $\eta_{\kappa,-2}$ to a holomorphic map $\tilde\eta_{\kappa,-2}:\pi_a^{-1}D_{-2}\to\pi_{a_\kappa}^{-1}D_{\kappa,-2}$.  We also lift $\eta_{\kappa,-1}=f_\kappa^{-1}\circ\eta_{\kappa,-2}\circ f$ to
$$	\tilde\eta_{\kappa,-1}=\tilde f_{\kappa,-2}^{-1}\circ\eta_{\kappa,-2}\circ\tilde f      $$
where the notation is that of Section 12, with obvious modification.  We write
$$ \tilde\eta_{\kappa0}=\tilde f_{\kappa,-1}^{-1}\circ\tilde\eta_{\kappa,-1}\circ\tilde f_{-1}$$
and
$$	\tilde\eta_{\kappa\beta}
	=(f_\kappa|V_{\kappa\beta})^{-1}\circ\tilde\eta_{\kappa\alpha}\circ f|V_\beta      $$
if order $\beta>0$ and $fV_\beta=V_\alpha$.  We have defined $\eta_{\kappa\alpha}$ above for $\alpha=-1,-2$, and we let $\eta_{\kappa\alpha}=\tilde\eta_{\kappa\alpha}$ when order $\alpha\ge0$.
\medskip
	We introduce a map $\eta_\kappa:{\cal A}_{\kappa1}\to{\cal A}_1$ by
$$	\eta_\kappa(\phi_{\kappa\alpha})
	=((\phi_{\kappa\alpha}\circ\tilde\eta_{\kappa\alpha}).\eta'_{\kappa\alpha})      $$
so that ${\cal L}_\kappa^\times=(\eta_\kappa,{\bf1}){\cal L}_\kappa(\eta_\kappa^{-1},{\bf1})$ acts on ${\cal A}$.  Using the decomposition
$$	{\cal L}_\kappa=\pmatrix{{\cal L}_{\kappa0}+{\cal L}_{\kappa1}&{\cal L}_{\kappa2}\cr	{\cal L}_{\kappa3}&{\cal L}_{\kappa4}\cr}      $$
as in Section 12, we define $L_\kappa^\times$ on ${\cal A}_1$ by
$$	L_\kappa^\times(\phi_\alpha)
	=\eta_\kappa({\cal L}_{\kappa0}+{\cal L}_{\kappa1})\eta_\kappa^{-1}(\phi_\alpha)
	+(\eta_\kappa^{-1}\phi_\alpha)_0(c_\kappa).\eta_\kappa{\cal L}_{\kappa2}
(|{1\over2}f''_\kappa(c_\kappa)\prod_{k=0}^{n-1}f'_\kappa(f_\kappa^kb_k)|^{-1/2})      $$
$$	={\cal L}_0(\phi_\alpha)+\eta_\kappa{\cal L}_{\kappa1}\eta_\kappa^{-1}(\phi_\kappa)
	+\eta'_{\kappa0}(c_\kappa)^{-1}\phi_0(c_\kappa).\eta_\kappa{\cal L}_{\kappa2}
(|{1\over2}f''_\kappa(c_\kappa)\prod_{k=0}^{n-1}f'_\kappa(f_\kappa^kb_k)|^{-1/2})      $$
$L_\kappa^\times$ is a compact perturbation of ${\cal L}_{\kappa0}$, and has therefore essential spectral radius $\le\gamma^{-1}$.  If $(\phi_\alpha)$ is a (generalized) eigenfunction of $L_\kappa^\times$ to the eigenvalue $\mu$, then 
$$	((\phi_\alpha),\eta_{\kappa0}(c_\kappa)^{-1}\phi_0(c_\kappa).
(|{1\over2}f''_\kappa(c_\kappa)\prod_{k=0}^{n-1}f'_\kappa(f_\kappa^kb_k)|^{-1/2}))      $$
is a (generalized) eigenfunction of ${\cal L}_\kappa^\times$ to the same eigenvalue $\mu$.  We have thus a multiplicity-preserving bijection of the eigenvalues $\mu$ of $L_\kappa^\times$ and ${\cal L}_\kappa^\times$ when $|\mu|>\max(\gamma^{-1},\delta\alpha^{1/2})$.  In particular, $1$ is a simple eigenvalue of $L_\kappa^\times$ for the values of $\kappa$ considered (a compact neighborhood $K$ of $0$).
\medskip
	The operator $L_\kappa^\times$ acting on ${\cal A}_1$ depends continuously on $\kappa$.  [This is because $\hat\phi_{\kappa\alpha}$, $\chi_{\kappa0}$, $\chi_{\kappa n\alpha}$ depend continuously on $\kappa$ (in particular, the $\chi_{\kappa n\alpha}$ for large $n$ are uniformly small).  Note however that ${\cal L}_\kappa^\times$ does not depend continuously on $\kappa$ because the continuity of $f'_\kappa(f_\kappa^kb_\kappa)$ is not uniform in $k$].  There is $\epsilon>0$ such that $L_\kappa^\times$ has no eigenvalue $\mu_\kappa$ with $|\mu_\kappa-1|<\epsilon$ except the simple eigenvalue $1$ [otherwise the continuity of $\kappa\to L_\kappa^\times$ would imply that $1$ has multiplicity $>1$ for some $\kappa$].  Therefore, the $1$-dimensional projection corresponding to the eigenvalue $1$ of $L_\kappa^\times$ depends continuously on $\kappa$, and so does the eigenvector $\Phi_\kappa^\times=(\eta_\kappa,1)\Phi_\kappa^0$ of ${\cal L}_\kappa^\times$, where $\Phi_\kappa^0$ denotes the eigenvector the the eigenvalue $1$ of ${\cal L}_\kappa$ normalized so that $w_\kappa\Phi_\kappa^0\ge0$ and $\int w_\kappa\Phi_\kappa^0=1$, with the obvious definition of $w_\kappa$ (involving the spikes $\psi_{\kappa n}$ associated with $f_\kappa$).
\medskip
	Note that a number of results have been obtained earlier on the continuous dependence of the a.c.i.m. $\rho$ on parameters.  I am indebted to Viviane Baladi for communicating the references [25], [27], [15], and also [26].
\vfill\eject
\noindent
{\bf References.}
\medskip
[1] A. Avila, M. Lyubich, and W. de Melo  "Regular or stochastic dynamics in real analytic families of unimodal maps."  Invent. Math. {\bf 154},451-550(2003).

[2] V. Baladi  {\it Positive transfer operators and decay of correlations.}  World Scientific, Singapore, 2000.

[3] V. Baladi  "On the susceptibility function on piecewise expanding interval maps."  Commun. Math. Phys. {\bf 275},839-859(2007).

[4] V. Baladi and G. Keller  "Zeta functions and transfer operators for piecewise monotone transformations."  Commun. Math. Phys. {\bf 127},459-479(1990).

[5] V. Baladi and D. Smania  "Linear response formula for piecewise expanding unimodal maps."  To appear.

[6] M. Benedicks and L. Carleson  "On iterations of $1-ax^2$ on $(-1,1).$  Ann. Math. {\bf 122},1-25(1985).

[7] M. Benedicks and L. Carleson  "The dynamics of the H\'enon map."  Ann. Math. {\bf 133},73-169(1991).

[8] M. Benedicks and L.-S. Young  "Absolutely continuous invariant measures and random perturbation for certain one-dimensional maps."  Ergod Th. Dynam Syst. {\bf 12},13-37(1992).

[9] O. Butterley and C. Liverani  "Smooth Anosov flows: correlation spectra and stability."  J. Modern Dynamics {\bf 1},301-322(2007).

[10] B. Cessac  "Does the complex susceptibility of the H\'enon map have a pole in the upper half plane?  A numerical investigation."  Submitted to Nonlinearity.

[11] L. Chierchia and G. Gallavotti  "Smooth prime integrals for quasi-integrable Hamiltonian systems."  Nuovo Cim. {\bf 67B},277-295(1982).

[12] P. Collet and J.-P. Eckmann  "Positive Lyapunov exponents and absolute continuity for maps of the interval."  Ergod Th. Dynam Syst. {\bf 3},13-46(1981).

[13] M. Jakobson  "Absolutely continuous invariant measures for certain maps of an interval."  Commun. Math. Phys. {\bf 81},39-88(1981).

[14] D. Dolgopyat  "On differentiability of SRB states for partially hyperbolic systems."  Invent. Math. {\bf 155},389-449(2004).

[15] J.M. Freitas  "Continuity of SRB measure and entropy for Benedicks-Carleson quad\-ratic maps."  Nonlinearity {\bf 18},831-854(2005).

[16] Y. Jiang and D. Ruelle  "Analyticity of the susceptibility function for unimodal Markovian maps of the interval."  Nonlinearity {\bf 18},2447-2453(2005).

[17] A. Katok, G. Knieper, M. Pollicott, and H. Weiss.  "Differentiability and analyticity of topological entropy for Anosov and geodesic flows."  Invent. Math. {\bf 98},581-597(1989).

[18] G. Keller and T. Nowicki  "Spectral theory, zeta functions and the distribution of periodic points for Collet-Eckmann maps."  Commun. Math. Phys. {\bf 149},31-69(1992).

[19] M. Misiurewicz  "Absolutely continuous measures for certain maps of an interval."  Publ. Math. IHES {\bf 53},17-52(1981).

[20] J. P\"oschel  "Integrability of Hamiltonian systems on Cantor sets."  Commun in Pure and Appl. Math. {\bf 35},653-696(1982).

[21] D. Ruelle  "Differentiation of SRB states.''  Commun. Math. Phys. {\bf 187},227-241(1997); ``Correction and complements.''  Commun. Math. Phys. {\bf 234},185-190(2003).  

[22] D. Ruelle  "Differentiation of SRB states for hyperbolic flows."  Ergod. Th. Dynam. Syst., to appear.

[23] D.Ruelle  "Differentiating the absolutely continuous invariant measure of an interval map f with respect to f."  Commun.Math. Phys. {\bf 258},445-453(2005).

[24] D. Ruelle  "Application of hyperbolic dynamics to physics: some problems and conjectures.''  Bull. Amer. Math. Soc. (N.S.) {\bf 41},275-278(2004).

[25] M. Rychlik and E. Sorets  "Regularity and other properties of absolutely continuous invariant measures for the quadratic family."  Commun. Math. Phys. {\bf 150},217-236(1992).

[26] B. Szewc  "The Perron-Frobenius operator in spaces of smooth functions on an interval."  Ergod. Th. Dynam. Syst. {\bf 4},613-643(1984).

[27] M. Tsujii  "On continuity of Bowen-Ruelle-Sinai measures in families of one dimensional
maps."  Commun. Math. Phys. {\bf 177},1-11(1996).

[28] Q.-D. Wang and L.-S. Young  "Nonuniformly expanding 1D maps."  Commun. Math. Phys. {\bf 264},255-282(2006).

[29] H. Whitney  "Analytic expansions of differentiable functions defined in closed sets.''  Trans. Amer. Math. Soc. {\bf 36},63-89(1934).

[30] L.-S. Young  "Decay of correlations for quadratic maps."  Commun. Math. Phys. {\bf 146},123-138(1992).

[31] L.-S. Young  "What are SRB measures, and which dynamical systems have them?"  J. Statistical Phys. {\bf 108},733-754(2002).

\end